\documentclass{amsart}

\usepackage{amssymb,amsmath,amsthm,latexsym,booktabs}

\theoremstyle{definition}
\newtheorem{remark}{Remark}[section]

\newcommand{\pp}{\phantom{-}}

\begin{document}

\title[Boolean matrices of order 3]
{Structure of the rational monoid algebra for Boolean matrices of order 3}

\author{Murray R. Bremner}

\address{Department of Mathematics and Statistics, University of Saskatchewan, Canada}

\email{bremner@math.usask.ca}

\urladdr{math.usask.ca/~bremner}



\thanks{Partially supported by a Discovery Grant from NSERC}

\begin{abstract}
We use computer algebra to study the 512-dimensional associative algebra $\mathbb{Q} \mathcal{B}_3$,
the rational monoid algebra of $3 \times 3$ Boolean matrices.
We obtain a basis for the radical in bijection with the 42 non-regular elements of $\mathcal{B}_3$.
The center of the 470-dimensional semisimple quotient has dimension 14;
we use a splitting algorithm to find a basis of orthogonal primitive idempotents.
We show that the semisimple quotient is the direct sum of simple two-sided ideals isomorphic to 
matrix algebras $M_d(\mathbb{Q})$ for $d = 1, 1, 1, 2, 3, 3, 3, 3, 6, 6, 7, 9, 9, 12$.
We construct the irreducible representations of $\mathcal{B}_3$ over $\mathbb{Q}$
by calculating the representation matrices for a minimal set of generators.
\end{abstract}

\maketitle

\allowdisplaybreaks


\section{Introduction}

We write $\mathcal{B}_n$ for the monoid of $n \times n$ Boolean matrices with the usual Boolean matrix product ($1+1=1$);
equivalently, $\mathcal{B}_n$ is the monoid of binary relations on $n$ elements with the relative product:
  \[
  R \circ S = \{ \, (i,j) \mid \text{for some $k$ we have $(i,k) \in R$ and $(k,j) \in S$} \, \}. 
  \]
Konieczny \cite{Konieczny} has recently a complete proof of Devadze's theorem \cite{Devadze} 
on minimal sets of generators for $\mathcal{B}_n$.

The representation theory of $\mathcal{B}_n$ over the field $\mathbb{Q}$ of rational numbers has received little attention.
One reason for this is Preston's theorem \cite{Preston}, which states that any finite group is a maximal subgroup 
of $\mathcal{B}_n$ for some $n$; it is therefore not realistic to expect a uniform structure theory for 
the rational monoid algebra $\mathbb{Q} \mathcal{B}_n$.
(The simplest non-trivial case $n = 2$ appears in \cite{BremnerElBachraoui}.)

This paper studies the irreducible rational representations of $\mathcal{B}_n$ for $n = 3$ by determining 
the structure of the rational monoid algebra $\mathbf{A} = \mathbb{Q} \mathcal{B}_3$.
We use the computer algebra system Maple, together with a constructive approach to the classical structure theory
of finite dimensional associative algebras; see \cite{Bremner}.
 
Section \ref{prelimsection} recalls some basic results about the structure of the monoid $\mathcal{B}_3$.

Section \ref{radicalsection} determines the radical $\mathbf{R} \subset \mathbf{A}$; we find that its dimension is 42
and we use the LLL algorithm \cite{BremnerLLL} for lattice basis reduction to find a natural basis of $\mathbf{R}$ 
in bijection with the non-regular elements of $\mathcal{B}_3$.

Section \ref{centersection} determines the structure constants for the semisimple quotient $\mathbf{S} = \mathbf{A}/\mathbf{R}$,
and a basis for the 14-dimensional center $\mathbf{C} \subset \mathbf{S}$; we then apply the splitting algorithm of 
Ivanyos and R\'onyai \cite{IvanyosRonyai} to determine a new basis for $\mathbf{C}$ consisting of orthogonal 
primitive idempotents.

Section \ref{decompsection} determines the decomposition of $\mathbf{S}$ into a direct sum of simple two-sided ideals; 
in particular, we find that $\mathbb{Q}$ is the splitting field of $\mathcal{B}_3$.
The dimensions of the irreducible representations are $d = 1$, 1, 1, 2, 3, 3, 3, 3, 6, 6, 7, 9, 9, 12.
We then find a minimal left ideal in each simple two-sided ideal, and
construct an isomorphism of each simple two-sided ideal with the matrix algebra
$M_d(\mathbb{Q})$.

Section \ref{repmatsection} calculates the representation matrices for the minimal set of five generators of 
$\mathcal{B}_3$ obtained from Devadze's theorem \cite{Devadze,Konieczny}; 
the matrix entries are integers, and only the four representations of dimensions $d \ge 7$ are faithful.


\section{Preliminaries} \label{prelimsection}

\subsection{Lex order} 

We order the elements $( b_{ij} )$ of $\mathcal{B}_3$ lexicographically, regarding them as binary numerals
$b_{11} b_{12} b_{13} b_{21} b_{22} b_{23} b_{31} b_{32} b_{33} $;
that is, we identify each element with an integer $n \in \mathcal{I} = \{ 1, \dots, 512 \}$ 
using the following bijection:
  \[
  \lambda\colon ( b_{ij} ) \longmapsto 1 + \sum_{i=1}^3 \sum_{j=1}^3 b_{ij} 2^{8-[3(i-1)+(j-1)]}.
  \]
We write $[n] = \lambda^{-1}(n)$ for the element corresponding to $n$.
The multiplication table of $\mathcal{B}_3$ is then given by the function 
$\mu\colon \mathcal{I} \times \mathcal{I} \to \mathcal{I}$
defined by $\mu( p, q ) = \lambda( [p] [q] )$.
If there is no loss of clarity, we will write $[n]$ simply as $n$.
  
\subsection{Regular elements}

An element $x \in \mathcal{B}_3$ is called regular if and only if $xyx = x$ for some $y \in \mathcal{B}_3$.
There are 42 non-regular elements in $\mathcal{B}_3$; see Table \ref{nonregular}.
These elements correspond to the following integers:
  \[
  \mathcal{NR}
  =
  \left\{
  \begin{array}{rrrrrrrrr} 
   95, & 111, & 116, & 118, & 158, & 172, & 175, & 182, & 207, 
  \\
  214, & 230, & 231, & 235, & 237, & 242, & 245, & 286, & 287, 
  \\
  300, & 308, & 335, & 340, & 343, & 347, & 349, & 356, & 370,
  \\
  371, & 396, & 398, & 406, & 410, & 413, & 420, & 426, & 427;
  \\
  239, & 246, & 351, & 372, & 414, & 428.
  \end{array}
  \right\}
  \]
The first 36 elements have four 0s; each of these has one row (say $i$) with two 0s 
and one column (say $j$) with two 0s, and the $(i,j)$ entry is 0.
The remaining 6 elements have three 0s: these are the complements of the permutation matrices.
We write $\mathcal{R} = \mathcal{I} \setminus \mathcal{NR}$ for the set of indices of the 470 regular elements.

\begin{table}
\[
\begin{array}{cccccc}
\left[ \begin{smallmatrix} 0 & 0 & 1 \\ 0 & 1 & 1 \\ 1 & 1 & 0 \end{smallmatrix} \right] & 
\left[ \begin{smallmatrix} 0 & 0 & 1 \\ 1 & 0 & 1 \\ 1 & 1 & 0 \end{smallmatrix} \right] & 
\left[ \begin{smallmatrix} 0 & 0 & 1 \\ 1 & 1 & 0 \\ 0 & 1 & 1 \end{smallmatrix} \right] & 
\left[ \begin{smallmatrix} 0 & 0 & 1 \\ 1 & 1 & 0 \\ 1 & 0 & 1 \end{smallmatrix} \right] & 
\left[ \begin{smallmatrix} 0 & 1 & 0 \\ 0 & 1 & 1 \\ 1 & 0 & 1 \end{smallmatrix} \right] & 
\left[ \begin{smallmatrix} 0 & 1 & 0 \\ 1 & 0 & 1 \\ 0 & 1 & 1 \end{smallmatrix} \right] \\[6pt] 
\left[ \begin{smallmatrix} 0 & 1 & 0 \\ 1 & 0 & 1 \\ 1 & 1 & 0 \end{smallmatrix} \right] & 
\left[ \begin{smallmatrix} 0 & 1 & 0 \\ 1 & 1 & 0 \\ 1 & 0 & 1 \end{smallmatrix} \right] & 
\left[ \begin{smallmatrix} 0 & 1 & 1 \\ 0 & 0 & 1 \\ 1 & 1 & 0 \end{smallmatrix} \right] & 
\left[ \begin{smallmatrix} 0 & 1 & 1 \\ 0 & 1 & 0 \\ 1 & 0 & 1 \end{smallmatrix} \right] & 
\left[ \begin{smallmatrix} 0 & 1 & 1 \\ 1 & 0 & 0 \\ 1 & 0 & 1 \end{smallmatrix} \right] & 
\left[ \begin{smallmatrix} 0 & 1 & 1 \\ 1 & 0 & 0 \\ 1 & 1 & 0 \end{smallmatrix} \right] \\[6pt]
\left[ \begin{smallmatrix} 0 & 1 & 1 \\ 1 & 0 & 1 \\ 0 & 1 & 0 \end{smallmatrix} \right] & 
\left[ \begin{smallmatrix} 0 & 1 & 1 \\ 1 & 0 & 1 \\ 1 & 0 & 0 \end{smallmatrix} \right] &
\left[ \begin{smallmatrix} 0 & 1 & 1 \\ 1 & 1 & 0 \\ 0 & 0 & 1 \end{smallmatrix} \right] & 
\left[ \begin{smallmatrix} 0 & 1 & 1 \\ 1 & 1 & 0 \\ 1 & 0 & 0 \end{smallmatrix} \right] & 
\left[ \begin{smallmatrix} 1 & 0 & 0 \\ 0 & 1 & 1 \\ 1 & 0 & 1 \end{smallmatrix} \right] & 
\left[ \begin{smallmatrix} 1 & 0 & 0 \\ 0 & 1 & 1 \\ 1 & 1 & 0 \end{smallmatrix} \right] \\[6pt] 
\left[ \begin{smallmatrix} 1 & 0 & 0 \\ 1 & 0 & 1 \\ 0 & 1 & 1 \end{smallmatrix} \right] &
\left[ \begin{smallmatrix} 1 & 0 & 0 \\ 1 & 1 & 0 \\ 0 & 1 & 1 \end{smallmatrix} \right] & 
\left[ \begin{smallmatrix} 1 & 0 & 1 \\ 0 & 0 & 1 \\ 1 & 1 & 0 \end{smallmatrix} \right] & 
\left[ \begin{smallmatrix} 1 & 0 & 1 \\ 0 & 1 & 0 \\ 0 & 1 & 1 \end{smallmatrix} \right] &
\left[ \begin{smallmatrix} 1 & 0 & 1 \\ 0 & 1 & 0 \\ 1 & 1 & 0 \end{smallmatrix} \right] & 
\left[ \begin{smallmatrix} 1 & 0 & 1 \\ 0 & 1 & 1 \\ 0 & 1 & 0 \end{smallmatrix} \right] \\[6pt]
\left[ \begin{smallmatrix} 1 & 0 & 1 \\ 0 & 1 & 1 \\ 1 & 0 & 0 \end{smallmatrix} \right] & 
\left[ \begin{smallmatrix} 1 & 0 & 1 \\ 1 & 0 & 0 \\ 0 & 1 & 1 \end{smallmatrix} \right] & 
\left[ \begin{smallmatrix} 1 & 0 & 1 \\ 1 & 1 & 0 \\ 0 & 0 & 1 \end{smallmatrix} \right] &
\left[ \begin{smallmatrix} 1 & 0 & 1 \\ 1 & 1 & 0 \\ 0 & 1 & 0 \end{smallmatrix} \right] & 
\left[ \begin{smallmatrix} 1 & 1 & 0 \\ 0 & 0 & 1 \\ 0 & 1 & 1 \end{smallmatrix} \right] & 
\left[ \begin{smallmatrix} 1 & 1 & 0 \\ 0 & 0 & 1 \\ 1 & 0 & 1 \end{smallmatrix} \right] \\[6pt]
\left[ \begin{smallmatrix} 1 & 1 & 0 \\ 0 & 1 & 0 \\ 1 & 0 & 1 \end{smallmatrix} \right] &
\left[ \begin{smallmatrix} 1 & 1 & 0 \\ 0 & 1 & 1 \\ 0 & 0 & 1 \end{smallmatrix} \right] &
\left[ \begin{smallmatrix} 1 & 1 & 0 \\ 0 & 1 & 1 \\ 1 & 0 & 0 \end{smallmatrix} \right] & 
\left[ \begin{smallmatrix} 1 & 1 & 0 \\ 1 & 0 & 0 \\ 0 & 1 & 1 \end{smallmatrix} \right] & 
\left[ \begin{smallmatrix} 1 & 1 & 0 \\ 1 & 0 & 1 \\ 0 & 0 & 1 \end{smallmatrix} \right] & 
\left[ \begin{smallmatrix} 1 & 1 & 0 \\ 1 & 0 & 1 \\ 0 & 1 & 0 \end{smallmatrix} \right] \\ \midrule
\left[ \begin{smallmatrix} 0 & 1 & 1 \\ 1 & 0 & 1 \\ 1 & 1 & 0 \end{smallmatrix} \right] & 
\left[ \begin{smallmatrix} 0 & 1 & 1 \\ 1 & 1 & 0 \\ 1 & 0 & 1 \end{smallmatrix} \right] &
\left[ \begin{smallmatrix} 1 & 0 & 1 \\ 0 & 1 & 1 \\ 1 & 1 & 0 \end{smallmatrix} \right] &
\left[ \begin{smallmatrix} 1 & 0 & 1 \\ 1 & 1 & 0 \\ 0 & 1 & 1 \end{smallmatrix} \right] & 
\left[ \begin{smallmatrix} 1 & 1 & 0 \\ 0 & 1 & 1 \\ 1 & 0 & 1 \end{smallmatrix} \right] & 
\left[ \begin{smallmatrix} 1 & 1 & 0 \\ 1 & 0 & 1 \\ 0 & 1 & 1 \end{smallmatrix} \right]
\end{array}
\]
\smallskip
\caption{Non-regular elements of $\mathcal{B}_3$}
\label{nonregular}
\end{table}

\subsection{$D$-classes}
 
Since $\mathcal{B}_3$ is finite, the $D$-classes coincide with the $J$-classes; 
thus $x, y \in \mathcal{B}_3$ belong to the same $D$-class
if and only if they generate the same two-sided monoid ideal, $\mathcal{B}_3 x \mathcal{B}_3 = \mathcal{B}_3 y \mathcal{B}_3$,
and this defines an equivalence relation on $\mathcal{B}_3$.
If a $D$-class contains a regular element then every element of that class is regular.
An element $x \in \mathcal{B}_3$ is called prime if it is not a permutation matrix and whenever $x = yz$ 
either $x$ or $y$ is a permutation matrix; see de Caen and Gregory \cite{deCaenGregory}.
If a $D$-class contains a prime element then every element of that class is prime.
The $D$-classes in $\mathcal{B}_3$ are displayed in Table \ref{dclasses}: 
each $D$-class is represented by its minimal element in lex order, together with its size,
and whether it is regular or prime.

\newcommand{\nnn}{\!\!\!\!\!\!\!}

\begin{table}
\begin{tabular}{cccccccccccc}
class &\!\!\!\!\!\!\!\!
$\left[ \begin{smallmatrix} 0 & 0 & 0 \\ 0 & 0 & 0 \\ 0 & 0 & 0 \end{smallmatrix} \right]$ &\nnn
$\left[ \begin{smallmatrix} 0 & 0 & 0 \\ 0 & 0 & 0 \\ 0 & 0 & 1 \end{smallmatrix} \right]$ &\nnn
$\left[ \begin{smallmatrix} 0 & 0 & 0 \\ 0 & 0 & 1 \\ 0 & 1 & 0 \end{smallmatrix} \right]$ &\nnn
$\left[ \begin{smallmatrix} 0 & 0 & 0 \\ 0 & 0 & 1 \\ 0 & 1 & 1 \end{smallmatrix} \right]$ &\nnn
$\left[ \begin{smallmatrix} 0 & 0 & 1 \\ 0 & 1 & 0 \\ 1 & 0 & 0 \end{smallmatrix} \right]$ &\nnn
$\left[ \begin{smallmatrix} 0 & 0 & 1 \\ 0 & 1 & 0 \\ 1 & 0 & 1 \end{smallmatrix} \right]$ &\nnn
$\left[ \begin{smallmatrix} 0 & 0 & 1 \\ 0 & 1 & 0 \\ 1 & 1 & 1 \end{smallmatrix} \right]$ &\nnn
$\left[ \begin{smallmatrix} 0 & 0 & 1 \\ 0 & 1 & 1 \\ 1 & 0 & 1 \end{smallmatrix} \right]$ &\nnn
$\left[ \begin{smallmatrix} 0 & 0 & 1 \\ 0 & 1 & 1 \\ 1 & 1 & 0 \end{smallmatrix} \right]$ &\nnn
$\left[ \begin{smallmatrix} 0 & 0 & 1 \\ 0 & 1 & 1 \\ 1 & 1 & 1 \end{smallmatrix} \right]$ &\nnn
$\left[ \begin{smallmatrix} 0 & 1 & 1 \\ 1 & 0 & 1 \\ 1 & 1 & 0 \end{smallmatrix} \right]$ 
\\[8pt]
size &
1 & 49 & 162 & 144 & 6 & 36 & 18 & 18 & 36 & 36 & 6
\\
regular &
yes & yes & yes & yes & yes & yes & yes & yes & no & yes & no
\\
prime &
no & no & no & no & no & no & no & no & no & no & yes
\end{tabular}
\smallskip
\caption{$D$-classes in $\mathcal{B}_3$}
\label{dclasses}
\end{table}

\subsection{Generators}

By the theorem of Devadze \cite{Devadze}, recently proved by Konieczny \cite{Konieczny},
the monoid $\mathcal{B}_3$ has the following minimal set of five generators:
  \begin{equation}
  \label{fivegenerators}
  \left[ \begin{smallmatrix} 0 & 1 & 0 \\ 1 & 0 & 0 \\ 0 & 0 & 1 \end{smallmatrix} \right], \qquad
  \left[ \begin{smallmatrix} 0 & 1 & 0 \\ 0 & 0 & 1 \\ 1 & 0 & 0 \end{smallmatrix} \right], \qquad
  \left[ \begin{smallmatrix} 1 & 0 & 0 \\ 1 & 1 & 0 \\ 0 & 0 & 1 \end{smallmatrix} \right], \qquad
  \left[ \begin{smallmatrix} 1 & 0 & 0 \\ 0 & 1 & 0 \\ 0 & 0 & 0 \end{smallmatrix} \right], \qquad
  \left[ \begin{smallmatrix} 0 & 1 & 1 \\ 1 & 0 & 1 \\ 1 & 1 & 0 \end{smallmatrix} \right].
  \end{equation}
The first two elements generate the subgroup of permutation matrices, the first four elements 
are regular but not prime, and the last element is prime but not regular.


\section{Non-regular elements and the radical} \label{radicalsection}

To find the radical of the monoid algebra $\mathbf{A} = \mathbb{Q} \mathcal{B}_3$, we use the theorem of Dickson \cite{Dickson}
together with Drazin's generalization of Maschke's theorem \cite{Drazin}.
Following \cite[Corollary 13]{Bremner}, we construct the $512 \times 512$ matrix $\Delta$ in which the $(i,j)$ entry is
  \[
  \Delta_{ij} = | \, \{ \, k \mid \mu( \, \mu( j, i ), \, k \, ) = k \, \} \, |.
  \]
We note that $\Delta$ is symmetric; its entries belong to $\{ 0, 1, 8, 27, 64, 125, 216, 512 \}$.  
The radical $\mathbf{R}$ of $\mathbf{A}$ consists of the linear combinations of the elements of $\mathcal{B}_3$ determined by the vectors
in the nullspace of $\Delta$.
Using the command \texttt{Rank} from the Maple package \texttt{LinearAlgebra} we find that $\Delta$ has rank 470 
and hence nullity 42.
Since the entries of $\Delta$ are integers, the most appropriate method to find a basis of its nullspace is 
by computing its Hermite normal form $H$ together with a unimodular integer matrix $U$ for which $U \Delta = H$; 
see \cite[Chapter 14]{BremnerLLL}.
Then the bottom 42 rows of $U$ form a lattice basis for the integer nullspace.
We use the following Maple command, which also applies the LLL algorithm to reduce the nullspace basis: 
  \[
  \texttt{HermiteForm( Delta, output=['U'], method='integer[reduced]' ):}
  \]
We obtain a basis for the radical in which all components belong to $\{ -1, 0, 1 \}$, but
the numbers of nonzero components in the vectors are 8 (17 times), 12 (14 times), 14 (twice), 16 (3 times),
18, 20 (twice), 22, 24 (twice).
Applying the LLL algorithm with a higher value of the parameter ($\alpha = 99/100$) 
produces an integer basis of the nullspace in which every vector has 8 nonzero components in $\{ \pm 1 \}$; 
see Table \ref{radicalbasis}.

\begin{table}
\footnotesize
\[
\begin{array}{rr}
 95 {-}  96 {-} 127 {+} 128 {-} 223 {+} 224 {+} 255 {-} 256, &\quad 
111 {-} 112 {-} 127 {+} 128 {-} 367 {+} 368 {+} 383 {-} 384, \\ 
116 {-} 120 {-} 124 {+} 128 {-} 244 {+} 248 {+} 252 {-} 256, &\quad 
118 {-} 120 {-} 126 {+} 128 {-} 374 {+} 376 {+} 382 {-} 384, \\ 
158 {-} 160 {-} 190 {+} 192 {-} 222 {+} 224 {+} 254 {-} 256, &\quad 
172 {-} 176 {-} 188 {+} 192 {-} 236 {+} 240 {+} 252 {-} 256, \\ 
175 {-} 176 {-} 191 {+} 192 {-} 431 {+} 432 {+} 447 {-} 448, &\quad
182 {-} 184 {-} 190 {+} 192 {-} 438 {+} 440 {+} 446 {-} 448, \\ 
207 {-} 208 {-} 223 {+} 224 {-} 463 {+} 464 {+} 479 {-} 480, &\quad 
214 {-} 216 {-} 222 {+} 224 {-} 470 {+} 472 {+} 478 {-} 480, \\ 
230 {-} 232 {-} 238 {+} 240 {-} 486 {+} 488 {+} 494 {-} 496, &\quad 
231 {-} 232 {-} 247 {+} 248 {-} 487 {+} 488 {+} 503 {-} 504, \\ 
235 {-} 236 {-} 251 {+} 252 {-} 491 {+} 492 {+} 507 {-} 508, &\quad 
237 {-} 238 {-} 253 {+} 254 {-} 493 {+} 494 {+} 509 {-} 510, \\ 
242 {-} 244 {-} 250 {+} 252 {-} 498 {+} 500 {+} 506 {-} 508, &\quad 
245 {-} 247 {-} 253 {+} 255 {-} 501 {+} 503 {+} 509 {-} 511, \\ 
286 {-} 288 {-} 318 {+} 320 {-} 350 {+} 352 {+} 382 {-} 384, &\quad 
287 {-} 288 {-} 319 {+} 320 {-} 415 {+} 416 {+} 447 {-} 448, \\ 
300 {-} 304 {-} 316 {+} 320 {-} 364 {+} 368 {+} 380 {-} 384, &\quad 
308 {-} 312 {-} 316 {+} 320 {-} 436 {+} 440 {+} 444 {-} 448, \\ 
335 {-} 336 {-} 367 {+} 368 {-} 463 {+} 464 {+} 495 {-} 496, &\quad 
340 {-} 344 {-} 348 {+} 352 {-} 468 {+} 472 {+} 476 {-} 480, \\ 
343 {-} 344 {-} 375 {+} 376 {-} 471 {+} 472 {+} 503 {-} 504, &\quad 
347 {-} 348 {-} 379 {+} 380 {-} 475 {+} 476 {+} 507 {-} 508, \\ 
349 {-} 350 {-} 381 {+} 382 {-} 477 {+} 478 {+} 509 {-} 510, &\quad 
356 {-} 360 {-} 364 {+} 368 {-} 484 {+} 488 {+} 492 {-} 496, \\
370 {-} 374 {-} 378 {+} 382 {-} 498 {+} 502 {+} 506 {-} 510, &\quad 
371 {-} 375 {-} 379 {+} 383 {-} 499 {+} 503 {+} 507 {-} 511, \\ 
396 {-} 400 {-} 412 {+} 416 {-} 460 {+} 464 {+} 476 {-} 480, &\quad 
398 {-} 400 {-} 430 {+} 432 {-} 462 {+} 464 {+} 494 {-} 496, \\ 
406 {-} 408 {-} 438 {+} 440 {-} 470 {+} 472 {+} 502 {-} 504, &\quad 
410 {-} 412 {-} 442 {+} 444 {-} 474 {+} 476 {+} 506 {-} 508, \\ 
413 {-} 415 {-} 445 {+} 447 {-} 477 {+} 479 {+} 509 {-} 511, &\quad 
420 {-} 424 {-} 436 {+} 440 {-} 484 {+} 488 {+} 500 {-} 504, \\
426 {-} 430 {-} 442 {+} 446 {-} 490 {+} 494 {+} 506 {-} 510, &\quad 
427 {-} 431 {-} 443 {+} 447 {-} 491 {+} 495 {+} 507 {-} 511, \\ \midrule 
239 {-} 240 {-} 255 {+} 256 {-} 495 {+} 496 {+} 511 {-} 512, &\quad
246 {-} 248 {-} 254 {+} 256 {-} 502 {+} 504 {+} 510 {-} 512, \\ 
351 {-} 352 {-} 383 {+} 384 {-} 479 {+} 480 {+} 511 {-} 512, &\quad 
372 {-} 376 {-} 380 {+} 384 {-} 500 {+} 504 {+} 508 {-} 512, \\ 
414 {-} 416 {-} 446 {+} 448 {-} 478 {+} 480 {+} 510 {-} 512, &\quad 
428 {-} 432 {-} 444 {+} 448 {-} 492 {+} 496 {+} 508 {-} 512.
\end{array}
\]
\smallskip 
\caption{Lattice reduced basis of the radical $\mathbf{R} \subset \mathbf{A}$}
\label{radicalbasis}
\end{table}

The leading terms of these basis elements are precisely the non-regular elements of $\mathcal{B}_3$.
These basis elements have a particular structure; we illustrate for one element of each type
(whether the leading term has three or four 0s).
The element with leading term $[95]$ in Table \ref{radicalbasis} represents this element of 
$\mathbb{Q} \mathcal{B}_3$:
\[
  \left[ \begin{smallmatrix} 0 & 0 & 1 \\ 0 & 1 & 1 \\ 1 & 1 & 0 \end{smallmatrix} \right]
- \left[ \begin{smallmatrix} 0 & 0 & 1 \\ 0 & 1 & 1 \\ 1 & 1 & 1 \end{smallmatrix} \right]
- \left[ \begin{smallmatrix} 0 & 0 & 1 \\ 1 & 1 & 1 \\ 1 & 1 & 0 \end{smallmatrix} \right]
+ \left[ \begin{smallmatrix} 0 & 0 & 1 \\ 1 & 1 & 1 \\ 1 & 1 & 1 \end{smallmatrix} \right]
- \left[ \begin{smallmatrix} 0 & 1 & 1 \\ 0 & 1 & 1 \\ 1 & 1 & 0 \end{smallmatrix} \right]
+ \left[ \begin{smallmatrix} 0 & 1 & 1 \\ 0 & 1 & 1 \\ 1 & 1 & 1 \end{smallmatrix} \right]
+ \left[ \begin{smallmatrix} 0 & 1 & 1 \\ 1 & 1 & 1 \\ 1 & 1 & 0 \end{smallmatrix} \right]
- \left[ \begin{smallmatrix} 0 & 1 & 1 \\ 1 & 1 & 1 \\ 1 & 1 & 1 \end{smallmatrix} \right]
\]
The other 7 terms are obtained by changing 0s to 1s,
except that the 0 in the upper left corner does not change; this is the only 0 which belongs to a row with two 0s
and a column with two 0s.
The sign of each term is $(-1)^i$ where $i$ is the number of 0s that have changed to 1s.
The element with leading term $[239]$ in Table \ref{radicalbasis} represents this element of 
$\mathbb{Q} \mathcal{B}_3$:
\[
  \left[ \begin{smallmatrix} 0 & 1 & 1 \\ 1 & 0 & 1 \\ 1 & 1 & 0 \end{smallmatrix} \right]
- \left[ \begin{smallmatrix} 0 & 1 & 1 \\ 1 & 0 & 1 \\ 1 & 1 & 1 \end{smallmatrix} \right]
- \left[ \begin{smallmatrix} 0 & 1 & 1 \\ 1 & 1 & 1 \\ 1 & 1 & 0 \end{smallmatrix} \right]
+ \left[ \begin{smallmatrix} 0 & 1 & 1 \\ 1 & 1 & 1 \\ 1 & 1 & 1 \end{smallmatrix} \right]
- \left[ \begin{smallmatrix} 1 & 1 & 1 \\ 1 & 0 & 1 \\ 1 & 1 & 0 \end{smallmatrix} \right]
+ \left[ \begin{smallmatrix} 1 & 1 & 1 \\ 1 & 0 & 1 \\ 1 & 1 & 1 \end{smallmatrix} \right]
+ \left[ \begin{smallmatrix} 1 & 1 & 1 \\ 1 & 1 & 1 \\ 1 & 1 & 0 \end{smallmatrix} \right]
- \left[ \begin{smallmatrix} 1 & 1 & 1 \\ 1 & 1 & 1 \\ 1 & 1 & 1 \end{smallmatrix} \right]
\]
As before, the other 7 terms are obtained by changing $i$ of the 0s to 1s, and the sign of each term is $(-1)^i$.


\section{The semisimple quotient and its center} \label{centersection}

\subsection{Structure constants}

Our description of the radical $\mathbf{R} \subset \mathbf{A}$ in terms of the non-regular elements of $\mathcal{B}_3$
allows us to construct easily the semisimple quotient $\mathbf{S} = \mathbf{A}/\mathbf{R}$.
As basis for $\mathbf{S}$ we take the (cosets modulo $\mathbf{R}$ of) the regular elements of $\mathcal{B}_3$;
that is, the elements $[p] + \mathbf{R}$ for $p \in \mathcal{R}$.
If $p, q \in \mathcal{R}$ the we have two cases.
If $[p][q]$ is also a regular element, then 
$\big( \, [p]+\mathbf{R} \, \big) \big( \, [q]+\mathbf{R} \, \big) = [p][q] + \mathbf{R}$.
If $[p][q]$ is not a regular element, then $[p][q] = [r_1]$ is the leading term of one of the radical basis elements 
  in Table \ref{radicalbasis}, which we write in the form 
  \[
  [r_1] - [r_2] - [r_3] + [r_4] - [r_5] + [r_6] + [r_7] - [r_8],
  \] 
  where $[r_2], \dots, [r_8]$ are regular elements, and therefore
   \[
  \big( \, [p]+\mathbf{R} \, \big) \big( \, [q]+\mathbf{R} \, \big) 
  = 
  \big( \, [r_2] + [r_3] - [r_4] + [r_5] - [r_6] - [r_7] + [r_8] \, \big) + \mathbf{R}.
  \]   
Let $s\colon \{ 1, \dots, 470 \} \to \mathcal{R}$ be the unique order preserving
bijection, so that a basis for $\mathbf{S}$ consists of the cosets
$\{ \, [s(k)] + \mathbf{R} \mid k = 1, \dots, 470 \, \}$.
We write the structure constants of $\mathbf{S}$ with respect to this basis as follows:
  \[
  \big( \, [ s(i) ] + \mathbf{R} \, \big) \big( \, [ s(j) ] + \mathbf{R} \, \big) 
  =
  \sum_{k=1}^{470} d_{ij}^n [ s(k) ] + \mathbf{R}. 
  \]   

\begin{table}
\footnotesize
\[
\begin{array}{rl}
%
A = &\!\!\!\!
   1, 
\\
%
%
B = &\!\!\!\! 
274, 
\\
%
%
C = &\!\!\!\! 
312 + 318 - 320 + 378 - 382 + 408 - 440 + 468 - 472 + 474 - 476 - 506 + 512, 
\\
%
%
D = &\!\!\!\! 
 54 -  56 -  62 +  64 - 310 + 345 - 346 - 377 + 382 + 388 - 392 - 404 + 440 - 452 + 456 
\\
&\!\!\!\!
{} - 473 + 476 + 505 - 512, 
\\
%
%
E = &\!\!\!\! 
210 - 212 - 218 + 220 + 279 - 280 + 298 - 302 - 311 - 314 + 320 - 362 + 366 - 407  
\\
&\!\!\!\!
{}+ 439 - 466 + 472 + 506 - 512, 
\\
%
%
F = &\!\!\!\! 
 47 -  48 -  63 +  64 + 199 - 200 - 215 + 216 + 233 - 234 - 249 + 250 - 303 + 304 + 319 
\\
&\!\!\!\!
{} - 320 - 455 + 456 + 471 - 472 - 489 + 490 + 505 - 506, 
\\
%
%
G = &\!\!\!\! 
109 - 110 - 125 + 126 + 167 - 168 - 183 + 184 + 203 - 204 - 219 + 220 - 365 + 366 
\\
&\!\!\!\!
{} + 381 - 382 - 423 + 424 + 439 - 440 - 459 + 460 + 475 - 476, 
\\
%
%
H = &\!\!\!\! 
  7 -   8 +  37 +  41 -  46 -  55 -  57 +  64 +  73 -  74 + 131 - 147 + 193 - 196 - 217 + 220 
\\
&\!\!\!\!
{} - 293 - 361 + 366 - 391 + 439 - 449 + 456 + 505 - 512, 
\\
%
%
I = &\!\!\!\! 
276 + 278 - 280 + 282 - 284 + 306 - 310 - 314 + 320 + 338 - 342 - 346 + 352 + 382  
\\
&\!\!\!\!
{}- 384 + 402 - 404 - 434 + 440 + 444 - 448 - 466 + 472 + 476 - 480 + 502 - 504 + 506 
\\
&\!\!\!\!
{} - 508 - 510 + 512, 
\\
%
%
J = &\!\!\!\! 
280 + 284 - 288 + 310 + 314 - 316 + 342 - 344 + 346 - 348 - 350 + 352 - 374 + 376 
\\
&\!\!\!\!
{} + 380 - 384 + 404 - 412 + 416 + 434 - 436 - 438 - 442 + 444 + 446 - 448 + 466 - 470 
\\
&\!\!\!\!
{} + 478 - 480 - 498 + 500 + 502 - 504 - 508 - 510 + 2 \cdot 512, 
\\
%
%
K = &\!\!\!\! 
 39 -  40 +  45 -  46 -  55 +  56 -  61 +  62 + 105 - 106 - 121 + 122 + 135 - 136 - 151 + 152 
\\
&\!\!\!\!
{} + 195 - 196 + 201 - 202 - 211 + 212 - 217 + 218 - 295 + 296 - 301 + 302 + 311 - 312 
\\
&\!\!\!\!
{} + 317 - 318 - 361 + 362 + 377 - 378 - 391 + 392 + 407 - 408 - 451 + 452 - 457 + 458 
\\
&\!\!\!\!
{} + 467 - 468 + 473 - 474, 
\\
%
%
L = &\!\!\!\! 
 79 -  80 + 125 - 126 - 127 + 128 + 168 + 171 - 176 - 184 - 187 + 192 + 215 - 216 - 223 
\\
&\!\!\!\!
{} + 224 + 229 - 232 + 234 - 236 - 238 + 240 - 247 + 248 - 250 + 252 - 253 + 254 + 255 
\\
&\!\!\!\!
{} - 256 + 303 - 304 - 319 + 320 - 367 + 368 - 381 + 382 + 383 - 384 - 424 - 431 + 432 
\\
&\!\!\!\!
{} + 440 + 447 - 448 + 459 - 460 - 463 + 464 - 471 + 472 - 475 + 476 + 479 - 480 - 485 
\\
&\!\!\!\!
{} + 488 - 490 - 491 + 492 + 494 + 495 - 496 + 503 - 504 + 506 + 507 - 508 + 509 - 510
\\
&\!\!\!\!
{} - 511, 
\\
%
%
M = &\!\!\!\! 
 99 - 100 - 103 + 104 - 107 + 108 - 115 + 119 + 123 - 128 + 141 - 142 - 143 + 144 - 157 
\\
&\!\!\!\!
{} + 159 - 173 + 174 + 189 - 192 - 205 + 206 + 221 - 224 - 227 + 228 + 243 - 248 - 252 
\\
&\!\!\!\!
{} - 254 + 2 \cdot 256 - 355 + 359 + 363 - 368 - 383 + 384 - 397 + 399 + 429 - 432 - 447 + 448 
\\
&\!\!\!\!
{} + 461 - 464 - 479 + 480 + 483 - 488 - 492 - 494 + 2 \cdot 496 - 503 + 504 - 507 + 508 - 509 
\\
&\!\!\!\!
{} + 510 + 2 \cdot 511 - 2 \cdot 512, 
\\
%
%
N = &\!\!\!\! 
 85 -  86 -  87 +  88 -  93 +  94 -  99 + 100 + 103 - 104 + 107 - 108 + 115 - 117 - 123 + 125 
\\
&\!\!\!\!
{} - 141 + 142 + 143 - 144 + 157 - 159 + 162 - 164 - 166 + 168 - 170 + 173 - 178 + 180 
\\
&\!\!\!\!
{} + 186 - 189 + 205 - 206 - 213 + 215 - 226 + 227 + 234 - 240 - 243 + 248 + 254 - 255 
\\
&\!\!\!\!
{} + 267 - 268 - 271 + 272 + 276 + 278 - 280 + 282 - 283 - 299 + 303 + 306 - 310 - 314 
\\
&\!\!\!\!
{} + 315 - 331 + 332 + 338- 341 - 346 + 352 + 355 - 359 + 373 - 376 - 380 + 383 - 395 
\\
&\!\!\!\!
{}
 + 397 + 402 - 404 + 411 - 416 - 418 + 422 - 429 + 432 + 444 - 446 + 459 - 461 - 466 
\\
&\!\!\!\!
{} + 469 - 478 + 479 + 482 - 483 + 492 - 495 - 500 + 502. 
\end{array}
\]
\caption{Lattice reduced basis of the center $\mathbf{C} \subset \mathbf{S}$}
\label{centerbasis}
\end{table}

\subsection{The center}

The next step is to compute the center $\mathbf{C}$ of the 470-dimensional semisimple associative algebra $\mathbf{S}$.
According to \cite[Corollary 15]{Bremner}, the coefficient vectors of the elements in $\mathbf{C}$ form the
nullspace of the $470^2 \times 470$ matrix $Z$ in which the entry in row $470(i{-}1) + k$ and column $j$ is
$d_{ij}^k - d_{ji}^k$.
Since $470^2 = 220900$ is so large, we do this calculation using modular arithmetic to save memory;
we use $p = 101$.
If $M$ is an integer matrix, then its rank over the finite field $\mathbb{F}_p$ is less than or equal to its rank 
over $\mathbb{Q}$, and so its nullspace over $\mathbb{F}_p$ contains (possibly strictly) the reduction modulo $p$
of its integer nullspace.
However, we will be able to verify that the basis we obtain for $\mathbf{C}$
using modular arithmetic spans a subalgebra of $\mathbf{C}$.
Thus we have upper and lower bounds for the dimension of $\mathbf{C}$ which are equal.

We first construct a modular matrix of size $470 \times 220900$, and store $d_{ij}^k - d_{ji}^k$ (mod $p$)
in row $j$ and column $470(i{-}1) + k$ for all $i, j, k = 1, \dots, 470$; this is $Z^t$ (mod $p$).
We use the command \texttt{RowReduce} from the Maple package \texttt{LinearAlgebra[Modular]} to compute 
the rank and the row canonical form (RCF).
The rank is 456 and so the nullity is 14; this is the dimension of the center.
We identify the 456 columns which contain the leading 1s of the rows of the RCF and record the corresponding pairs 
$(i_1,k_1), \dots, (i_{456},k_{456})$.
In this way we have identified the rows of $Z$ which form a basis for the row space.

We next construct an integer matrix $Z_0$ of size $456 \times 470$, and store 
$d_{i_\ell j}^{k_\ell} - d_{j i_\ell}^{k_\ell}$ in row $\ell$ and column $j$;
by the previous paragraph, the reductions mod $p$ of the rows form a basis 
of the row space of $Z$ mod $p$.
We use the Maple command
  \[
  \texttt{HermiteForm( Transpose(Z0), output='U', method='integer[reduced]' ):}
  \]
and extract the bottom 14 rows forming a lattice reduced basis of the center $\mathbf{C}$.
(Using a higher value of the parameter in the LLL algorithm does not improve the results in this case.)
We sort these rows by increasing Euclidean length, and change signs where necessary to ensure that
each row has a positive leading entry.
The components of the vectors in this basis of $\mathbf{C}$ belong to $\{ 0, \pm 1, \pm 2 \}$,
and the squared lengths of the vectors are 1, 1, 13, 19, 19, 24, 24, 25, 31, 40, 48, 72, 72, 95.
These basis vectors $A, \dots, N$ are displayed in Table \ref{centerbasis}, and
the multiplication table for this ordered basis appears in Table \ref{centerstructure}.

\newcommand{\nn}{\!\!\!\!\!\!}

\begin{table}
\footnotesize
\[
\begin{array}{cccccccccccccc} 
A &\nn A &\nn A &\nn {-}A &\nn {-}A &\nn 0 &\nn 0 &\nn {-}A &\nn A &\nn 2A &\nn 0 &\nn 0 &\nn 0 &\nn A  \\ 
A &\nn B &\nn C &\nn D &\nn E &\nn F &\nn G &\nn H &\nn I &\nn J &\nn K &\nn L &\nn M &\nn N  \\ 
A &\nn C &\nn C &\nn {-}C &\nn {-}C &\nn 0 &\nn 0 &\nn H &\nn C &\nn 2C &\nn K &\nn 0 &\nn 0 &\nn C  \\ 
{-}A &\nn D &\nn {-}C &\nn {-}D &\nn C &\nn {-}F &\nn 0 &\nn {-}H &\nn D &\nn {-}C{+}D &\nn {-}K &\nn F &\nn 0 &\nn D  \\ 
{-}A &\nn E &\nn {-}C &\nn C &\nn {-}E &\nn 0 &\nn {-}G &\nn {-}H &\nn E &\nn {-}C{+}E &\nn {-}K &\nn G &\nn 0 &\nn E  \\ 
0 &\nn F &\nn 0 &\nn {-}F &\nn 0 &\nn {-}C{-}D &\nn 0 &\nn 0 &\nn F &\nn F &\nn 0 &\nn C{+}D &\nn 0 &\nn F  \\ 
0 &\nn G &\nn 0 &\nn 0 &\nn {-}G &\nn 0 &\nn {-}C{-}E &\nn 0 &\nn G &\nn G &\nn 0 &\nn C{+}E &\nn 0 &\nn G  \\ 
{-}A &\nn H &\nn H &\nn {-}H &\nn {-}H &\nn 0 &\nn 0 &\nn {-}H &\nn H &\nn 2H &\nn 0 &\nn 0 &\nn 0 &\nn H  \\ 
A &\nn I &\nn C &\nn D &\nn E &\nn F &\nn G &\nn H &\nn I &\nn J &\nn K &\nn L &\nn 0 &\nn I  \\ 
2A &\nn J &\nn 2C &\nn {-}C{+}D &\nn {-}C{+}E &\nn F &\nn G &\nn 2H &\nn J &\nn 2C{+}J &\nn 2K &\nn L &\nn 0 &\nn J  \\ 
0 &\nn K &\nn K &\nn {-}K &\nn {-}K &\nn 0 &\nn 0 &\nn 0 &\nn K &\nn 2K &\nn {-}K &\nn 0 &\nn 0 &\nn K  \\ 
0 &\nn L &\nn 0 &\nn F &\nn G &\nn C{+}D &\nn C{+}E &\nn 0 &\nn L &\nn L &\nn 0 &\nn {-}2C{+}J &\nn 0 &\nn L  \\ 
0 &\nn M &\nn 0 &\nn 0 &\nn 0 &\nn 0 &\nn 0 &\nn 0 &\nn 0 &\nn 0 &\nn 0 &\nn 0 &\nn 2B{-}2I{+}M &\nn {-}2B{+}M{+}2N  \\ 
A &\nn N &\nn C &\nn D &\nn E &\nn F &\nn G &\nn H &\nn I &\nn J &\nn K &\nn L &\nn {-}2B{+}M{+}2N &\nn 5B{-}4N
\end{array}
\]
\smallskip
\caption{Multiplication table for the center $\mathbf{C} \subset \mathbf{S}$}
\label{centerstructure}
\end{table}

\subsection{Orthogonal primitive idempotents}
 
Our next task is to find a basis for $\mathbf{C}$ consisting of orthogonal primitive idempotents;
these idempotents play the role of the identity matrices in the simple two-sided ideals $M_d(\mathbb{Q})$
in the Wedderburn decomposition of $\mathbf{S}$.
We use the splitting algorithm of Ivanyos and R\'onyai \cite{IvanyosRonyai}; see also 
\cite[\S 1.8]{Bremner}.
This algorithm repeatedly splits a central ideal $X$ into the orthogonal direct sum of two central ideals $Y$ and $Z$;
thus, $X = Y \oplus Z$ and $YZ = \{0\}$.
The new basis will be denoted by the first 14 letters of the Greek alphabet: $\alpha, \beta, \dots, \xi$.
Throughout this computation it is important to remember that scalar terms in polynomials must be interpreted
as scalar multiples of the identity element in $X$.

\subsubsection{Splitting 1.}
We start with $X$ equal to the center $\mathbf{C}$ with identity element $B$.
We let $x = A$ (or any element that is not a scalar multiply of the identity).
Using Table \ref{centerstructure} we find that the minimal polynomial of $x$ is $x(x-1)$.
We therefore set $y = x = A$ and $z = x-1 = A-B$ and find that the ideal $Y = \langle y \rangle$ is 1-dimensional 
and the ideal ideal $Z = \langle z \rangle$ is 13-dimensional;
moreover, $X = Y \oplus Z$ and $YZ = \{0\}$.
Since $A^2 = A$, we record $\alpha = A$ as a primitive idempotent.

\subsubsection{Splitting 2.}
We start with $X = \langle A-B \rangle$; this is the ideal $Z = \langle z \rangle$ from Splitting 1
with identity element $-A+B$.
We compute a basis for $X$ to find an element which is not a scalar multiple of the identity; we take $x = A-N$
with minimal polynomial $(x+1)(x-5)$.
We therefore set $y = x+1 = B-N$ and $z = x-5 = 6A-5B-N$; we find that $Y = \langle y \rangle$ is 1-dimensional
and $Z = \langle z \rangle$ is 12-dimensional.
We have $(B-N)^2 = 6(B-N)$ so we record $\beta = \tfrac16(B-N)$ as a primitive idempotent.

\subsubsection{Splitting 3.}
We start with $X = \langle 6A-5B-N \rangle$; the identity element is $-A + \tfrac56 B + \tfrac16 N$.
An element of $X$ which is not a scalar multiple of the identity is $x = A - \tfrac52 M - N$.
Its minimal polynomial is $(x+1) (x-\tfrac32) (x+6)$ so we set $y = x+1 = \tfrac56 B - \tfrac52 M - \tfrac56 N$
and $z = x^2 + \tfrac92 x - 9 = \tfrac{25}{2} A - \tfrac{25}{2} I$.
We find that $y$ generates a 2-dimensional ideal $Y$ and $z$ generates a 10-dimensional ideal $Z$.
We do not obtain a primitive idempotent; we first split $Y$ in Splitting 4 to obtain two primitive
idempotents, and then we split $Z$ in Splitting 5.

\subsubsection{Splitting 4.}
We start with $X = \langle \tfrac56 B - \tfrac52 M - \tfrac56 N \rangle$;
the identity element is $\tfrac56 B - I + \tfrac16 N$.
We take $x = B - 3 M - 3 N$; its minimal polynomial is $(x-3)(x+6)$.
We set $y = x-3 = -\tfrac32 B + 3 I - 3 M - \tfrac32 N$ and 
$z = x + 6 = 6 B - 6 I - 3 M$.
Both $y$ and $z$ generate 1-dimensional ideals; scaling gives two primitive idempotents,
$\gamma = \tfrac16 B - \tfrac13 I + \tfrac13 M + \tfrac16 N$ and $\delta = \tfrac23 B - \tfrac23 I - \tfrac13 M$.

\subsubsection{Splitting 5.}
We start with $X = \langle \tfrac{25}{2} A - \tfrac{25}{2} I \rangle$, the ideal $Z$ from Splitting 3.
The identity element is $-A-I$;
we choose $x = A - \tfrac12 J$ with minimal polynomial $x(x+\tfrac12)(x+1)$.
We choose $y = x = A - \tfrac12 J$ and 
$z = x^2 + \tfrac32 x + \tfrac12 = \tfrac12 C + \tfrac12 I - \tfrac12 J$.
We find that $Y = \langle y \rangle$ is 9-dimensional and $Z = \langle z \rangle$ is 1-dimensional;
scaling $z$ gives the primitive idempotent $\epsilon = C + I - J$.

\subsubsection{Splitting 6.}
We start with $X = \langle A - \tfrac12 J \rangle$, the ideal $Y$ from Splitting 5.
The identity element is $-A-C+J$; we choose $x = A - \tfrac12 J$, with minimal polynomial $(x+1)(x+\tfrac12)$.
We set $y = x+1 = -C+\tfrac12 J$ and $z = x+\tfrac12 = \tfrac12 A - \tfrac12 C$.
We find that $Y = \langle y \rangle$ is 6-dimensional and $Z = \langle z \rangle$ is 3-dimensional.
We first split $Z$ in Splittings 7 and 8, and then split $Y$ in Splitting 9.

\subsubsection{Splitting 7.}
We start with $X = \langle \tfrac12 A - \tfrac12 C \rangle$, the ideal $Z$ from Splitting 6.
The identity element is $-A+C$; we choose $x = A+H$, with minimal polynomial $x(x+1)$.
We set $y = x = A+H$ and $z = x+1 = C+H$.
The ideal $Y = \langle y \rangle$ is 1-dimensional, giving the new primitive idempotent
$\zeta = -A-H$; the ideal $Z = \langle z \rangle$ is 2-dimensional.

\subsubsection{Splitting 8.}
We start with $X = \langle C+H \rangle$.
The identity element is $C+H$, so we take $x = K$ with minimal polynomial $x(x+1)$.
We set $y = x = K$ and $z = x + 1 = C+H+K$, and find that both generate 1-dimensional ideals,
giving the new primitive idempotents $\eta = -K$ and $\theta = C+H+K$.

\subsubsection{Splitting 9.}
We start with $X = \langle -C+\tfrac12 J \rangle$, the ideal $Y$ from Splitting 6.
The identity element is $-2C + J$ so we take $x = D + \tfrac12 J$ with minimal polynomial $(x+\tfrac12)(x-\tfrac12)$.
We set $y = x+\tfrac12 = -C+D+J$ and $z = x-\tfrac12 = C+D$ and find that 
$Y = \langle y \rangle$ is 4-dimensional and $Z = \langle z \rangle$ is 2-dimensional.
 We first split $Z$ in Splitting 10, and then split $Y$ in Splitting 11.

\subsubsection{Splitting 10.}
We start with $X = \langle C+D \rangle$, the ideal $Z$ from Splitting 9, with identity element $-C-D$.
We take $x = F$ with minimal polynomial $(x+1)(x-1)$.
We set $y = x+1 = -C-D+F$ and $z = x-1 = C+D+F$; both generate 1-dimensional ideals and produce the new primtive
idempotents $\iota = -\tfrac12 C - \tfrac12 D + \tfrac12 F$ and $\kappa = -\tfrac12 C - \tfrac12 D + \tfrac12 F$.

\subsubsection{Splitting 11.}
We start with $X = \langle -C+D+J \rangle$, the ideal $Y$ from Splitting 9, with identity element $-C+D+J$.
We take $x = C+E$ with minimal polynomial $x(x+1)$.
We set $y = x = C+E$ and $z = x+1 = D+E+J$; these elements generate 2-dimensional ideals $Y$ and $Z$.

\subsubsection{Splitting 12.}
We start with $X = \langle C+E \rangle$, the ideal $Y$ from Splitting 11, with identity element $-C-E$.
We choose $x = G$, with minimal polynomial $(x+1)(x-1)$; we set $y = x+1 = -C-E+G$ and $z = x-1 = C+E+G$.
Both $y$ and $z$ generate 1-dimensional ideals, giving two new primitive idempotents
$\lambda = -\tfrac12 C - \tfrac12 E + \tfrac12 G$ and $\mu = -\tfrac12 C - \tfrac12 E - \tfrac12 G$.

\subsubsection{Splitting 13.}
We start with $X = \langle D+E+J \rangle$, the ideal $Z$ from Splitting 11, with identity element $D+E+J$.
We take $x = F+G+L$ with minimal polynomial $(x+1)(x-1)$.
We set $y = x+1 = D+E+F+G+J+L$ and $z = x-1 = -D-E+F+G-J+L$.
Both $y$ and $z$ generate 1-dimensional ideals, giving our final primitive idempotents
$\nu = \tfrac12 D+\tfrac12 E+\tfrac12 F+\tfrac12 G+\tfrac12 J+\tfrac12 L$ and 
$\xi = \tfrac12 D+\tfrac12 E-\tfrac12 F-\tfrac12 G+\tfrac12 J-\tfrac12 L$.

\subsubsection{Summary.}

We have calculated a new basis of the center $\mathbf{C} \subset \mathbf{S}$ consisting of the orthogonal
primitive idempotents whose coefficient vectors with respect to the basis of 
Table \ref{centerbasis} are given by the rows of the matrix in Table \ref{opitable}.
We check using Table \ref{centerstructure} that $x^2 = x$ and $xy = 0$ for all $x, y$ in this new basis. 
We note that all the minimal polynomials obtained in this calculation split into linear factors over 
$\mathbb{Q}$. 

\begin{table}
\footnotesize
\[
\frac16
\left[
\begin{array}{rrrrrrrrrrrrrr}  
   6 &  . &  . &  . &  . &  . &  . &  . &  . &  . &  . &  . &  . &  . \\
   . &  1 &  . &  . &  . &  . &  . &  . &  . &  . &  . &  . &  . & -1 \\
   . &  1 &  . &  . &  . &  . &  . &  . & -2 &  . &  . &  . &  2 &  1 \\
   . &  4 &  . &  . &  . &  . &  . &  . & -4 &  . &  . &  . & -2 &  . \\
   . &  . &  6 &  . &  . &  . &  . &  . &  6 & -6 &  . &  . &  . &  . \\
  -6 &  . &  . &  . &  . &  . &  . & -6 &  . &  . &  . &  . &  . &  . \\
   . &  . &  . &  . &  . &  . &  . &  . &  . &  . & -6 &  . &  . &  . \\
   . &  . &  6 &  . &  . &  . &  . &  6 &  . &  . &  6 &  . &  . &  . \\
   . &  . & -3 & -3 &  . &  3 &  . &  . &  . &  . &  . &  . &  . &  . \\
   . &  . & -3 & -3 &  . & -3 &  . &  . &  . &  . &  . &  . &  . &  . \\
   . &  . & -3 &  . & -3 &  . &  3 &  . &  . &  . &  . &  . &  . &  . \\
   . &  . & -3 &  . & -3 &  . & -3 &  . &  . &  . &  . &  . &  . &  . \\
   . &  . &  . &  3 &  3 &  3 &  3 &  . &  . &  3 &  . &  3 &  . &  . \\
   . &  . &  . &  3 &  3 & -3 & -3 &  . &  . &  3 &  . & -3 &  . &  .
\end{array}
\right]
\]
\smallskip
\caption{Orthogonal primitive idempotents in the center $\mathbf{C} \subset \mathbf{S}$}
\label{opitable}
\end{table}


\section{The decomposition into simple ideals} \label{decompsection}

\begin{table}
\scriptsize
\[
\begin{array}{rl}
%
A = &\!\!\!\!
    1, 
\\[2pt]
%
%
B = &\!\!\!\!
\tfrac16
(
{-}  85 {+}  86 {+}  87 {-}  88 {+}  93 {-}  94 {+}  99 {-} 100 {-} 103 {+} 104 
{-} 107 {+} 108 {-} 115 {+} 117 {+} 123 {-} 125 {+} 141 {-} 142 
\\
&\!\!\!\! {}
{-} 143 {+} 144 
{-} 157 {+} 159 {-} 162 {+} 164 {+} 166 {-} 168 {+} 170 {-} 173 {+} 178 {-} 180 
{-} 186 {+} 189 {-} 205 {+} 206 {+} 213 
\\
&\!\!\!\! {}
{-} 215 {+} 226 {-} 227 {-} 234 {+} 240 
{+} 243 {-} 248 {-} 254 {+} 255 {-} 267 {+} 268 {+} 271 {-} 272 {+} 274 {-} 276 
{-} 278 {+} 280 
\\
&\!\!\!\! {}
{-} 282 {+} 283 {+} 299 {-} 303 {-} 306 {+} 310 {+} 314 {-} 315 
{+} 331 {-} 332 {-} 338 {+} 341 {+} 346 {-} 352 {-} 355 {+} 359 {-} 373 
\\
&\!\!\!\! {}
{+} 376 
{+} 380 {-} 383 {+} 395 {-} 397 {-} 402 {+} 404 {-} 411 {+} 416 {+} 418 {-} 422 
{+} 429 {-} 432 {-} 444 {+} 446 {-} 459 {+} 461 
\\
&\!\!\!\! {}
{+} 466 {-} 469 {+} 478 {-} 479 
{-} 482 {+} 483 {-} 492 {+} 495 {+} 500 {-} 502
), 
\\[2pt]
%
%
C = &\!\!\!\!
\tfrac16
(
   85 {-}  86 {-}  87 {+}  88 {-}  93 {+}  94 {+}  99 {-} 100 {-} 103 {+} 104 
{-} 107 {+} 108 {-} 115 {-} 117 {+} 2 \cdot 119 {+} 123 {+} 125 
\\
&\!\!\!\! {}
{-} 2 \cdot 128 {+} 141 {-} 142 
{-} 143 {+} 144 {-} 157 {+} 159 {+} 162 {-} 164 {-} 166 {+} 168 {-} 170 {-} 173 
{+} 2 \cdot 174 {-} 178 {+} 180 
\\
&\!\!\!\! {}
{+} 186 {+} 189 {-} 2 \cdot 192 {-} 205 {+} 206 {-} 213 {+} 215 
{+} 2 \cdot 221 {-} 2 \cdot 224 {-} 226 {-} 227 {+} 2 \cdot 228 {+} 234 {-} 240 {+} 243 
\\
&\!\!\!\! {}
{-} 248 {-} 2 \cdot 252 
{-} 254 {-} 255 {+} 4 \cdot 256 {+} 267 {-} 268 {-} 271 {+} 272 {+} 274 {-} 276 {-} 278 
{+} 280 {-} 282 {-} 283 
\\
&\!\!\!\! {}
{+} 2 \cdot 284 {-} 299 {+} 303 {-} 306 {+} 310 {+} 314 {+} 315 
{-} 2 \cdot 320 {-} 331 {+} 332 {-} 338 {-} 341 {+} 2 \cdot 342 {+} 346 {-} 352 
\\
&\!\!\!\! {}
{-} 355 {+} 359 
{+} 2 \cdot 363 {-} 2 \cdot 368 {+} 373 {-} 376 {-} 380 {-} 2 \cdot 382 {-} 383 {+} 4 \cdot 384 {-} 395 {-} 397 
{+} 2 \cdot 399 {-} 402 
\\
&\!\!\!\! {}
{+} 404 {+} 411 {-} 416 {-} 418 {+} 422 {+} 429 {-} 432 {+} 2 \cdot 434 
{-} 2 \cdot 440 {-} 444 {-} 446 {-} 2 \cdot 447 {+} 4 \cdot 448 {+} 459 {+} 461 
\\
&\!\!\!\! {}
{-} 2 \cdot 464 {+} 466 {+} 469 
{-} 2 \cdot 472 {-} 2 \cdot 476 {-} 478 {-} 479 {+} 4 \cdot 480 {+} 482 {+} 483 {-} 2 \cdot 488 {-} 492 {-} 2 \cdot 494 
{-} 495 
\\
&\!\!\!\! {}
{+} 4 \cdot 496 {-} 500 {-} 502 {-} 2 \cdot 503 {+} 4 \cdot 504 {-} 2 \cdot 506 {-} 2 \cdot 507 {+} 4 \cdot 508 {-} 2 \cdot 509 
{+} 4 \cdot 510 {+} 4 \cdot 511 
\\
&\!\!\!\! {}
{-}6 \cdot 512
), 
\\[2pt]
%
%
D = &\!\!\!\!
\tfrac13
(
{-}  99 {+} 100 {+} 103 {-} 104 {+} 107 {-} 108 {+} 115 {-} 119 {-} 123 {+} 128 
{-} 141 {+} 142 {+} 143 {-} 144 {+} 157 {-} 159 
\\
&\!\!\!\! {}
{+} 173 {-} 174 {-} 189 {+} 192 
{+} 205 {-} 206 {-} 221 {+} 224 {+} 227 {-} 228 {-} 243 {+} 248 {+} 252 {+} 254 
{-} 2 \cdot 256 {+} 2 \cdot 274 
\\
&\!\!\!\! {}
{-} 2 \cdot 276 {-} 2 \cdot 278 {+} 2 \cdot 280 {-} 2 \cdot 282 {+} 2 \cdot 284 {-} 2 \cdot 306 {+} 2 \cdot 310 {+} 2 \cdot 314 
{-} 2 \cdot 320 {-} 2 \cdot 338 {+} 2 \cdot 342 
\\
&\!\!\!\! {}
{+} 2 \cdot 346 {-} 2 \cdot 352 {+} 355 {-} 359 {-} 363 {+} 368 {-} 2 \cdot 382 
{+} 383 {+} 384 {+} 397 {-} 399 {-} 2 \cdot 402 {+} 2 \cdot 404 {-} 429 
\\
&\!\!\!\! {}
{+} 432 {+} 2 \cdot 434 {-} 2 \cdot 440 
{-} 2 \cdot 444 {+} 447 {+} 448 {-} 461 {+} 464 {+} 2 \cdot 466 {-} 2 \cdot 472 {-} 2 \cdot 476 {+} 479 {+} 480 
\\
&\!\!\!\! {}
{-} 483 {+} 488 {+} 492 {+} 494 {-} 2 \cdot 496 {-} 2 \cdot 502 {+} 503 {+} 504 {-} 2 \cdot 506 {+} 507 
{+} 508 {+} 509 {+} 510 {-} 2 \cdot 511
), 
\\[2pt]
%
%
K = &\!\!\!\!
\tfrac12
(
  109 {-} 110 {-} 125 {+} 126 {+} 167 {-} 168 {-} 183 {+} 184 {+} 203 {-} 204 
{-} 210 {+} 212 {+} 218 {-} 219 {-} 279 {+} 280 
\\
&\!\!\!\! {}
{-} 298 {+} 302 {+} 311 {-} 312 
{+} 314 {-} 318 {+} 362 {-} 365 {-} 378 {+} 381 {+} 407 {-} 408 {-} 423 {+} 424 
{-} 459 {+} 460 {+} 466 
\\
&\!\!\!\! {}
{-} 468 {-} 474 {+} 475
), 
\\[2pt]
%
%
I = &\!\!\!\!
\tfrac12
(
   47 {-}  48 {-}  54 {+}  56 {+}  62 {-}  63 {+} 199 {-} 200 {-} 215 {+} 216 
{+} 233 {-} 234 {-} 249 {+} 250 {-} 303 {+} 304 {+} 310 {-} 312 
\\
&\!\!\!\! {}
{-} 318 {+} 319 
{-} 345 {+} 346 {+} 377 {-} 378 {-} 388 {+} 392 {+} 404 {-} 408 {+} 452 {-} 455 
{-} 468 {+} 471 {+} 473 {-} 474 {-} 489 
\\
&\!\!\!\! {}
{+} 490
), 
\\[2pt]
%
%
L = &\!\!\!\!
\tfrac12
(
{-} 109 {+} 110 {+} 125 {-} 126 {-} 167 {+} 168 {+} 183 {-} 184 {-} 203 {+} 204 
{-} 210 {+} 212 {+} 218 {+} 219 {-} 2 \cdot 220 {-} 279 
\\
&\!\!\!\! {}
{+} 280 {-} 298 {+} 302 {+} 311 
{-} 312 {+} 314 {-} 318 {+} 362 {+} 365 {-} 2 \cdot 366 {-} 378 {-} 381 {+} 2 \cdot 382 {+} 407 
{-} 408 {+} 423 
\\
&\!\!\!\! {}
{-} 424 {-} 2 \cdot 439 {+} 2 \cdot 440 {+} 459 {-} 460 {+} 466 {-} 468 {-} 474 
{-} 475 {+} 2 \cdot 476
), 
\\[2pt]
%
%
J = &\!\!\!\!
\tfrac12
(
{-}  47 {+}  48 {-}  54 {+}  56 {+}  62 {+}  63 {-} 2 \cdot  64 {-} 199 {+} 200 {+} 215 
{-} 216 {-} 233 {+} 234 {+} 249 {-} 250 {+} 303 {-} 304 
\\
&\!\!\!\! {}
{+} 310 {-} 312 {-} 318 
{-} 319 {+} 2 \cdot 320 {-} 345 {+} 346 {+} 377 {-} 378 {-} 388 {+} 392 {+} 404 {-} 408 
{+} 452 {+} 455 {-} 2 \cdot 456 
\\
&\!\!\!\! {}
{-} 468 {-} 471 {+} 2 \cdot 472 {+} 473 {-} 474 {+} 489 {-} 490 
{-} 2 \cdot 505 {+} 2 \cdot 506
), 
\\[2pt]
%
%
G = &\!\!\!\!
{-}  39 {+}  40 {-}  45 {+}  46 {+}  55 {-}  56 {+}  61 {-}  62 {-} 105 {+} 106 
{+} 121 {-} 122 {-} 135 {+} 136 {+} 151 {-} 152 {-} 195 {+} 196 
\\
&\!\!\!\! {}
{-} 201 {+} 202 
{+} 211 {-} 212 {+} 217 {-} 218 {+} 295 {-} 296 {+} 301 {-} 302 {-} 311 {+} 312 
{-} 317 {+} 318 {+} 361 {-} 362 {-} 377 
\\
&\!\!\!\! {}
{+} 378 {+} 391 {-} 392 {-} 407 {+} 408 
{+} 451 {-} 452 {+} 457 {-} 458 {-} 467 {+} 468 {-} 473 {+} 474, 
\\[2pt]
%
%
E = &\!\!\!\!
  276 {+} 278 {-} 2 \cdot 280 {+} 282 {-} 2 \cdot 284 {+} 288 {+} 306 {-} 2 \cdot 310 {+} 312 {-} 2 \cdot 314 
{+} 316 {+} 318 {+} 338 {-} 2 \cdot 342 
\\
&\!\!\!\! {}
{+} 344 {-} 2 \cdot 346 {+} 348 {+} 350 {+} 374 {-} 376 
{+} 378 {-} 380 {+} 402 {-} 2 \cdot 404 {+} 408 {+} 412 {-} 416 {-} 2 \cdot 434 {+} 436 
\\
&\!\!\!\! {}
{+} 438 
{+} 442 {-} 446 {-} 2 \cdot 466 {+} 468 {+} 470 {+} 474 {-} 478 {+} 498 {-} 500, 
\\[2pt]
%
%
F = &\!\!\!\!
{-}   1 {-}   7 {+}   8 {-}  37 {-}  41 {+}  46 {+}  55 {+}  57 {-}  64 {-}  73 
{+}  74 {-} 131 {+} 147 {-} 193 {+} 196 {+} 217 {-} 220 {+} 293 {+} 361 {-} 366 
\\
&\!\!\!\! {}
{+} 391 {-} 439 {+} 449 {-} 456 {-} 505 {+} 512, 
\\[2pt]
%
%
N = &\!\!\!\!
\tfrac12
(
{-}  47 {+}  48 {+}  54 {-}  56 {-}  62 {+}  63 {-}  79 {+}  80 {-} 109 {+} 110 
{+} 127 {-} 128 {-} 167 {-} 171 {+} 176 {+} 183 {+} 187 {-} 192 
\\
&\!\!\!\! {}
{-} 199 {+} 200 
{-} 203 {+} 204 {+} 210 {-} 212 {-} 218 {+} 219 {+} 223 {-} 224 {-} 229 {+} 232 
{-} 233 {+} 236 {+} 238 {-} 240 {+} 247 
\\
&\!\!\!\! {}
{-} 248 {+} 249 {-} 252 {+} 253 {-} 254 
{-} 255 {+} 256 {+} 279 {+} 284 {-} 288 {+} 298 {-} 302 {-} 311 {-} 316 {+} 320 
{+} 342 {-} 344 
\\
&\!\!\!\! {}
{+} 345 {-} 348 {-} 350 {+} 352 {-} 362 {+} 365 {+} 367 {-} 368 
{-} 374 {+} 376 {-} 377 {+} 380 {+} 382 {-} 383 {+} 388 {-} 392 {-} 407 
\\
&\!\!\!\! {}
{-} 412 
{+} 416 {+} 423 {+} 431 {-} 432 {+} 434 {-} 436 {-} 438 {+} 440 {-} 442 {+} 444 
{+} 446 {-} 447 {-} 452 {+} 455 {+} 463 {-} 464 
\\
&\!\!\!\! {}
{-} 470 {+} 472 {-} 473 {+} 476 
{+} 478 {-} 479 {+} 485 {-} 488 {+} 489 {+} 491 {-} 492 {-} 494 {-} 495 {+} 496 
{-} 498 {+} 500 {+} 502 
\\
&\!\!\!\! {}
{-} 503 {+} 506 {-} 507 {-} 509 {+} 511
), 
\\[2pt]
%
%
M = &\!\!\!\!
\tfrac12
(
   47 {-}  48 {+}  54 {-}  56 {-}  62 {-}  63 {+} 2 \cdot  64 {+}  79 {-}  80 {+} 109 
{-} 110 {-} 127 {+} 128 {+} 167 {+} 171 {-} 176 {-} 183 {-} 187 
\\
&\!\!\!\! {}
{+} 192 {+} 199 
{-} 200 {+} 203 {-} 204 {+} 210 {-} 212 {-} 218 {-} 219 {+} 2 \cdot 220 {-} 223 {+} 224 
{+} 229 {-} 232 {+} 233 {-} 236 
\\
&\!\!\!\! {}
{-} 238 {+} 240 {-} 247 {+} 248 {-} 249 {+} 252 
{-} 253 {+} 254 {+} 255 {-} 256 {+} 279 {+} 284 {-} 288 {+} 298 {-} 302 {-} 311 
{-} 316 
\\
&\!\!\!\! {}
{+} 320 {+} 342 {-} 344 {+} 345 {-} 348 {-} 350 {+} 352 {-} 362 {-} 365 
{+} 2 \cdot 366 {-} 367 {+} 368 {-} 374 {+} 376 {-} 377 {+} 380 
\\
&\!\!\!\! {}
{+} 382 {+} 383 {-} 2 \cdot 384 
{+} 388 {-} 392 {-} 407 {-} 412 {+} 416 {-} 423 {-} 431 {+} 432 {+} 434 {-} 436 
{-} 438 {+} 2 \cdot 439 {+} 440 
\\
&\!\!\!\! {}
{-} 442 {+} 444 {+} 446 {+} 447 {-} 2 \cdot 448 {-} 452 {-} 455 
{+} 2 \cdot 456 {-} 463 {+} 464 {-} 470 {+} 472 {-} 473 {+} 476 {+} 478 {+} 479 
\\
&\!\!\!\! {}
{-} 2 \cdot 480 
{-} 485 {+} 488 {-} 489 {-} 491 {+} 492 {+} 494 {+} 495 {-} 496 {-} 498 {+} 500 
{+} 502 {+} 503 {-} 2 \cdot 504 {+} 2 \cdot 505 
\\
&\!\!\!\! {}
{+} 506 {+} 507 {-} 2 \cdot 508 {+} 509 {-} 2 \cdot 510 {-} 511
), 
\\[2pt]
%
%
H = &\!\!\!\!
    7 {-}   8 {+}  37 {+}  39 {-}  40 {+}  41 {+}  45 {-} 2 \cdot  46 {-} 2 \cdot  55 {+}  56 
{-}  57 {-}  61 {+}  62 {+}  64 {+}  73 {-}  74 {+} 105 {-} 106 {-} 121 {+} 122 
\\
&\!\!\!\! {}
{+} 131 {+} 135 {-} 136 {-} 147 {-} 151 {+} 152 {+} 193 {+} 195 {-} 2 \cdot 196 {+} 201 
{-} 202 {-} 211 {+} 212 {-} 2 \cdot 217 {+} 218 {+} 220 
\\
&\!\!\!\! {}
{-} 293 {-} 295 {+} 296 {-} 301 
{+} 302 {+} 311 {+} 317 {-} 320 {-} 2 \cdot 361 {+} 362 {+} 366 {+} 377 {-} 382 {-} 2 \cdot 391 
{+} 392 {+} 407 
\\
&\!\!\!\! {}
{+} 439 {-} 440 {-} 449 {-} 451 {+} 452 {+} 456 {-} 457 {+} 458 
{+} 467 {-} 472 {+} 473 {-} 476 {+} 505 {-} 506. 
\end{array}
\]
\smallskip
\caption{Orthogonal idempotents (identity matrices) in $\mathbf{S}$}
\label{identitytable}
\end{table}

Each of the orthogonal primitive idempotents $\alpha, \dots, \xi$ in the rows of the matrix 
in Table \ref{opitable} represents a linear combination of the basis elements $A, \dots, N$ 
of the center $\mathbf{C} \subset \mathbf{S}$ in Table \ref{centerbasis}.
The full expansions of these orthogonal idempotents as linear combinations of the basis elements of $\mathbf{S}$
(that is, $[p] + R$ for $p \in \mathcal{R}$) are given in Table \ref{identitytable}.
These elements represent the identity matrices in the simple two-sided ideals $M_d(\mathbb{Q})$
in the Wedderburn decomposition of $\mathbf{S}$; since $\mathbf{C} \subset \mathbf{S}$ we use the same notation 
$A, \dots, N$ for these elements regarded as elements of $\mathbf{C}$ or $\mathbf{S}$.
However, the order of the basis elements has been changed in Table \ref{identitytable};
we now have $A, B, C, D, K, I, L, J, G, E, F, N, M, H$ for the following reason.

We calculate the dimension of the two-sided ideal of $\mathbf{S}$ generated by each element;
since these elements are central, it suffices to consider all left multiples by basis elements of $\mathbf{S}$.
We obtain the values 1, 1, 1, 4, 9, 9, 9, 9, 36, 36, 49, 81, 81, 144 for the elements in the order of 
Table \ref{identitytable}.
We have reordered the basis elements of the center so that these dimensions are non-decreasing.
We therefore expect that 
  \[
  \mathbf{S} \cong
  3 \, \mathbb{Q} \oplus
  M_2(\mathbb{Q}) \oplus
  4 \, M_3(\mathbb{Q}) \oplus
  2 \, M_6(\mathbb{Q}) \oplus
  M_7(\mathbb{Q}) \oplus
  2 \, M_9(\mathbb{Q}) \oplus
  M_{12}(\mathbb{Q}).
  \]
To verify this, we need to construct an isomorphism of each two-sided ideal with the appropriate matrix algebra,
and this requires finding a minimal left ideal inside each two-sided ideal.
For a general finite dimensional associative algebra over $\mathbb{Q}$, this problem is very difficult; 
see Ivanyos et al. \cite{IRS} for an algorithm which calls oracles for 
factoring integers and polynomials over finite fields.

On the other hand, in practice it is often easy to find an element of the two-sided ideal which generates
a minimal left ideal.
In the present case we are very lucky: for every two-sided ideal except one, the first basis vector in a 
row-reduced basis of the two-sided ideal generates a minimal left ideal.
The exceptional case is the two-sided ideal of dimension 4: none of the basis vectors of the two-sided ideal
generates a minimal left ideal, but the sum of the first two basis vectors does.

\begin{table}
\footnotesize
\[
\begin{array}{rl}
d \,\,\,\, &
\\
%
%
2\colon &
   85 -  86 -  87 +  88 -  93 +  94 +  99 - 100 - 103 + 104 - 107 + 108 - 115 - 117 + 2 \cdot 119 
\\
& {}
+ 123 + 125 - 2 \cdot 128 - 213 + 215 + 221 - 224 - 227 + 228 + 243 - 248 - 252 - 255 
\\
& {}
+ 2 \cdot 256 - 267 + 268 + 271 - 272 - 274 + 276 + 278 - 280 + 282 + 283 - 2 \cdot 284 + 299 
\\
& {}
- 303 + 306 - 310 - 314 - 315 + 2 \cdot 320 + 331 - 332 + 338 - 341 - 346 + 352 - 355 
\\
& {}
+ 359 + 373 - 376 + 380 - 383 + 395 - 399 + 402 - 404 - 411 + 416 - 434 + 440 + 444 
\\
& {}
+ 447 - 2 \cdot 448 - 459 + 464 - 466 + 469 + 2 \cdot 476 - 478 - 480 + 483 - 488 - 492 + 495 
\\
& {}
+ 502 - 2 \cdot 503 + 504 + 506 - 508 - 509 + 511,
\\
%
%
3\colon &
   84 -  88 -  92 +  96 - 140 + 144 + 156 - 160 + 204 - 208 - 212 + 216,
\\ 
%
%
3\colon &
   30 -  32 -  44 +  48 +  60 -  62 -  94 +  96 + 108 - 112 - 124 + 126,
\\ 
%
%
3\colon &
   84 -  88 -  92 +  96 + 140 - 144 - 156 + 160 - 204 + 208 - 212 + 216 + 2 \cdot 220 - 2 \cdot 224,
\\ 
%
%
3\colon &
   30 -  32 +  44 -  48 -  60 -  62 + 2 \cdot  64 -  94 +  96 - 108 + 112 + 124 + 126 - 2 \cdot 128,
\\ 
%
%
6\colon &
   12 -  16 -  28 +  32 -  76 +  80 +  92 -  96,
\\ 
%
%
6\colon &
   86 -  88 -  94 +  96 - 342 + 344 + 350 - 352,
\\ 
%
%
7\colon &
    1 - 512,
\\ 
%
%
9\colon & 
   11 -  12 -  18 +  20 +  26 -  27,
\\ 
%
%
9\colon & 
   11 -  12 +  18 -  20 -  26 -  27 + 2 \cdot  28,
\\ 
%
%
12\colon &
    2 -   8 -  74 +  80. 
\end{array}
\]
\caption{Generators of minimal left ideals for $d > 1$}
\label{minimal}
\end{table}

The first three two-sided ideals are 1-dimensional, so there is nothing to do: these ideals are isomorphic to 
(the $1 \times 1$ matrices over) $\mathbb{Q}$.
For the remaining 11 two-sided ideals, the generators of minimal left ideals are given in Table \ref{minimal}.
These elements are also generators of the two-sided ideal to which they belong.
The last six elements are especially simple, and their terms are of combinatorial interest, so we write
them out in full as rational linear combinations of Boolean matrices:
  \begin{align*}
  6\colon &\quad
    \left[ \begin{smallmatrix} 0&0&0 \\ 0&0&1 \\ 0&1&1 \end{smallmatrix} \right]
  - \left[ \begin{smallmatrix} 0&0&0 \\ 0&0&1 \\ 1&1&1 \end{smallmatrix} \right]
  - \left[ \begin{smallmatrix} 0&0&0 \\ 0&1&1 \\ 0&1&1 \end{smallmatrix} \right]
  + \left[ \begin{smallmatrix} 0&0&0 \\ 0&1&1 \\ 1&1&1 \end{smallmatrix} \right]
  - \left[ \begin{smallmatrix} 0&0&1 \\ 0&0&1 \\ 0&1&1 \end{smallmatrix} \right]
  + \left[ \begin{smallmatrix} 0&0&1 \\ 0&0&1 \\ 1&1&1 \end{smallmatrix} \right]
  + \left[ \begin{smallmatrix} 0&0&1 \\ 0&1&1 \\ 0&1&1 \end{smallmatrix} \right]
  - \left[ \begin{smallmatrix} 0&0&1 \\ 0&1&1 \\ 1&1&1 \end{smallmatrix} \right]
  \\
  6\colon &\quad
    \left[ \begin{smallmatrix} 0&0&1 \\ 0&1&0 \\ 1&0&1 \end{smallmatrix} \right]
  - \left[ \begin{smallmatrix} 0&0&1 \\ 0&1&0 \\ 1&1&1 \end{smallmatrix} \right]
  - \left[ \begin{smallmatrix} 0&0&1 \\ 0&1&1 \\ 1&0&1 \end{smallmatrix} \right]
  + \left[ \begin{smallmatrix} 0&0&1 \\ 0&1&1 \\ 1&1&1 \end{smallmatrix} \right]
  - \left[ \begin{smallmatrix} 1&0&1 \\ 0&1&0 \\ 1&0&1 \end{smallmatrix} \right]
  + \left[ \begin{smallmatrix} 1&0&1 \\ 0&1&0 \\ 1&1&1 \end{smallmatrix} \right]
  + \left[ \begin{smallmatrix} 1&0&1 \\ 0&1&1 \\ 1&0&1 \end{smallmatrix} \right] 
  - \left[ \begin{smallmatrix} 1&0&1 \\ 0&1&1 \\ 1&1&1 \end{smallmatrix} \right]
  \\
  7\colon &\quad
    \left[ \begin{smallmatrix} 0&0&0 \\ 0&0&0 \\ 0&0&0 \end{smallmatrix} \right] 
  - \left[ \begin{smallmatrix} 1&1&1 \\ 1&1&1 \\ 1&1&1 \end{smallmatrix} \right]
  \\
  9\colon &\quad
    \left[ \begin{smallmatrix} 0&0&0 \\ 0&0&1 \\ 0&1&0 \end{smallmatrix} \right]
  - \left[ \begin{smallmatrix} 0&0&0 \\ 0&0&1 \\ 0&1&1 \end{smallmatrix} \right]
  - \left[ \begin{smallmatrix} 0&0&0 \\ 0&1&0 \\ 0&0&1 \end{smallmatrix} \right]
  + \left[ \begin{smallmatrix} 0&0&0 \\ 0&1&0 \\ 0&1&1 \end{smallmatrix} \right]
  + \left[ \begin{smallmatrix} 0&0&0 \\ 0&1&1 \\ 0&0&1 \end{smallmatrix} \right] 
  - \left[ \begin{smallmatrix} 0&0&0 \\ 0&1&1 \\ 0&1&0 \end{smallmatrix} \right]
  \\
  9\colon &\quad
    \left[ \begin{smallmatrix} 0&0&0 \\ 0&0&1 \\ 0&1&0 \end{smallmatrix} \right]
  - \left[ \begin{smallmatrix} 0&0&0 \\ 0&0&1 \\ 0&1&1 \end{smallmatrix} \right]
  + \left[ \begin{smallmatrix} 0&0&0 \\ 0&1&0 \\ 0&0&1 \end{smallmatrix} \right]
  - \left[ \begin{smallmatrix} 0&0&0 \\ 0&1&0 \\ 0&1&1 \end{smallmatrix} \right]
  - \left[ \begin{smallmatrix} 0&0&0 \\ 0&1&1 \\ 0&0&1 \end{smallmatrix} \right]
  - \left[ \begin{smallmatrix} 0&0&0 \\ 0&1&1 \\ 0&1&0 \end{smallmatrix} \right]
  +2 \left[ \begin{smallmatrix} 0&0&0 \\ 0&1&1 \\ 0&1&1 \end{smallmatrix} \right]
  \\
  12\colon &\quad
    \left[ \begin{smallmatrix} 0&0&0 \\ 0&0&0 \\ 0&0&1 \end{smallmatrix} \right]
  - \left[ \begin{smallmatrix} 0&0&0 \\ 0&0&0 \\ 1&1&1 \end{smallmatrix} \right]
  - \left[ \begin{smallmatrix} 0&0&1 \\ 0&0&1 \\ 0&0&1 \end{smallmatrix} \right]
  + \left[ \begin{smallmatrix} 0&0&1 \\ 0&0&1 \\ 1&1&1 \end{smallmatrix} \right]
  \end{align*}
It is now straightforward to construct an isomorphism of each two-sided ideal $\mathbf{I}$ with the corresponding 
matrix algebra $M_d(\mathbb{Q})$.
Once we have a $d$-dimensional minimal left ideal $\mathbf{L} \subset \mathbf{I}$, 
we choose a basis for $\mathbf{L}$ and identify
the basis elements with the standard basis vectors $e_1, \dots, e_d \in \mathbb{Q}^d$.
Since $\mathbf{I}\,\mathbf{L} \subset \mathbf{L}$, we calculate the $d \times d$ matrix representing the left 
action on $\mathbf{L}$ of each of the $d^2$ basis elements of $\mathbf{I}$.
For each $i, j = 1, \dots, d$ we then solve for the matrix unit $E_{ij}$ as a linear combination of these
$d \times d$ matrices.
Equivalently, we set $E_{ij}$ to the general linear combination of the basis elements of $\mathbf{I}$
with indeterminate coefficients, and solve
the linear system determined by the equations $E_{ij} e_k = \delta_{jk} e_i$ for all $k = 1, \dots, d$.


\section{Representation matrices for the generators} \label{repmatsection}

Once we have calculated the elements of each $d^2$-dimensional simple two-sided ideal $\mathbf{I}$ representing 
the matrix units $E_{ij}$ ($i, j = 1, \dots, d$), it is straightforward to compute the matrix representing any
element $[p] \in \mathcal{B}_3$ ($p \in \mathcal{I}$) in each of the irreducible representations.
We first multiply the coset $[p] + \mathbf{R} \in \mathbf{S}$ by the element $I_d \in \mathbf{S}$ 
from Table \ref{identitytable} representing the identity matrix in the given two-sided ideal.
Then $R_d(p) = I_d\,( [p] + \mathbf{R} ) \in \mathbf{I} \cong M_d(\mathbb{Q})$, and we solve a linear system 
to express $R_d(p)$ as a linear combination of the elements of $\mathbf{I}$ representing the matrix units $E_{ij}$
($i, j = 1, \dots, d$).
This linear combination of matrix units is the matrix of $p$ in the irreducible representation of $\mathcal{B}_3$ 
corresponding to $\mathbf{I}$.
We perform this calculation for the generators \eqref{fivegenerators} of $\mathcal{B}_3$ displayed in 
Section \ref{prelimsection}, and obtain the following results.

\subsection{Representation 1; dimension 1}

We find that $R(p) = 1$ for all five generators, and hence for all $[p] \in \mathcal{B}_3$: 
the unit representation of $\mathcal{B}_3$.

\subsection{Representations 2, 3, 4; dimensions 1, 1, 2}

These are the sign, unit, and 2-dimensional irreducible representations of the subgroup 
$\mathcal{S}_3 \subset \mathcal{B}_3$; all other elements of $\mathcal{B}_3$ are sent to 0:
  \begin{alignat*}{4}
  &2\colon &\quad
  &\left[ \begin{smallmatrix} 0 & 1 & 0 \\ 1 & 0 & 0 \\ 0 & 0 & 1 \end{smallmatrix} \right] \mapsto -1 &\qquad
  &\left[ \begin{smallmatrix} 0 & 1 & 0 \\ 0 & 0 & 1 \\ 1 & 0 & 0 \end{smallmatrix} \right] \mapsto  1 &\qquad
  &\left[ \begin{smallmatrix} 1 & 0 & 0 \\ 1 & 1 & 0 \\ 0 & 0 & 1 \end{smallmatrix} \right]
   \left[ \begin{smallmatrix} 1 & 0 & 0 \\ 0 & 1 & 0 \\ 0 & 0 & 0 \end{smallmatrix} \right]
   \left[ \begin{smallmatrix} 0 & 1 & 1 \\ 1 & 0 & 1 \\ 1 & 1 & 0 \end{smallmatrix} \right] \mapsto  0
  \\
  &3\colon &\quad
  &\left[ \begin{smallmatrix} 0 & 1 & 0 \\ 1 & 0 & 0 \\ 0 & 0 & 1 \end{smallmatrix} \right] \mapsto  1 &\qquad
  &\left[ \begin{smallmatrix} 0 & 1 & 0 \\ 0 & 0 & 1 \\ 1 & 0 & 0 \end{smallmatrix} \right] \mapsto  1 &\qquad
  &\left[ \begin{smallmatrix} 1 & 0 & 0 \\ 1 & 1 & 0 \\ 0 & 0 & 1 \end{smallmatrix} \right]
   \left[ \begin{smallmatrix} 1 & 0 & 0 \\ 0 & 1 & 0 \\ 0 & 0 & 0 \end{smallmatrix} \right]
   \left[ \begin{smallmatrix} 0 & 1 & 1 \\ 1 & 0 & 1 \\ 1 & 1 & 0 \end{smallmatrix} \right] \mapsto  0
  \\
  &4\colon &\quad
  &\left[ \begin{smallmatrix} 0 & 1 & 0 \\ 1 & 0 & 0 \\ 0 & 0 & 1 \end{smallmatrix} \right] 
  \mapsto 
  \left[ \begin{smallmatrix} 1 & -1 \\ 0 & -1 \end{smallmatrix} \right] 
  &\qquad
  &\left[ \begin{smallmatrix} 0 & 1 & 0 \\ 0 & 0 & 1 \\ 1 & 0 & 0 \end{smallmatrix} \right] 
  \mapsto
  \left[ \begin{smallmatrix} 0 & -1 \\ 1 & -1 \end{smallmatrix} \right] 
  &\qquad
  &\left[ \begin{smallmatrix} 1 & 0 & 0 \\ 1 & 1 & 0 \\ 0 & 0 & 1 \end{smallmatrix} \right]
  \left[ \begin{smallmatrix} 1 & 0 & 0 \\ 0 & 1 & 0 \\ 0 & 0 & 0 \end{smallmatrix} \right]
  \left[ \begin{smallmatrix} 0 & 1 & 1 \\ 1 & 0 & 1 \\ 1 & 1 & 0 \end{smallmatrix} \right]
  \mapsto
  \left[ \begin{smallmatrix} 0 & 0 \\ 0 & 0 \end{smallmatrix} \right]
  \end{alignat*}

\subsection{Representations 5, 6, 7, 8; dimension 3}

In all four cases the last two generators are sent to the $3 \times 3$ zero matrix; 
the first three generators have the following representations:
  \begin{alignat*}{4}
  &
  5\colon &\quad
  &
  \left[ \begin{smallmatrix} 0 & 1 & 0 \\ 1 & 0 & 0 \\ 0 & 0 & 1 \end{smallmatrix} \right] 
  \mapsto 
  \left[ \begin{smallmatrix} -1 & 0 & 0 \\ \pp 0 & 0 & 1 \\ \pp 0 & 1 & 0 \end{smallmatrix} \right] 
  &\qquad
  &
  \left[ \begin{smallmatrix} 0 & 1 & 0 \\ 0 & 0 & 1 \\ 1 & 0 & 0 \end{smallmatrix} \right] 
  \mapsto
  \left[ \begin{smallmatrix} \pp 0 & \pp 0 & 1 \\ -1 & \pp 0 & 0 \\ \pp 0 & -1 & 0 \end{smallmatrix} \right] 
  &\qquad
  &
  \left[ \begin{smallmatrix} 1 & 0 & 0 \\ 1 & 1 & 0 \\ 0 & 0 & 1 \end{smallmatrix} \right]
  \mapsto
  \left[ \begin{smallmatrix} 0 & 0 & 0 \\ 0 & 1 & 0 \\ 0 & 0 & 0 \end{smallmatrix} \right]
  \\ 
  &
  6\colon &\quad
  &
  \left[ \begin{smallmatrix} 0 & 1 & 0 \\ 1 & 0 & 0 \\ 0 & 0 & 1 \end{smallmatrix} \right] 
  \mapsto 
  \left[ \begin{smallmatrix} 0 & 1 & \pp 0 \\ 1 & 0 & \pp 0 \\ 0 & 0 & -1 \end{smallmatrix} \right] 
  &\qquad
  &
  \left[ \begin{smallmatrix} 0 & 1 & 0 \\ 0 & 0 & 1 \\ 1 & 0 & 0 \end{smallmatrix} \right] 
  \mapsto
  \left[ \begin{smallmatrix} 0 & -1 & \pp 0 \\ 0 & \pp 0 & -1 \\ 1 & \pp 0 & \pp 0 \end{smallmatrix} \right] 
  &\qquad
  &
  \left[ \begin{smallmatrix} 1 & 0 & 0 \\ 1 & 1 & 0 \\ 0 & 0 & 1 \end{smallmatrix} \right]
  \mapsto
  \left[ \begin{smallmatrix} 1 & 0 & 0 \\ 0 & 0 & 0 \\ 0 & 0 & 0 \end{smallmatrix} \right]
  \\ 
  &
  7\colon &\quad
  &
  \left[ \begin{smallmatrix} 0 & 1 & 0 \\ 1 & 0 & 0 \\ 0 & 0 & 1 \end{smallmatrix} \right] 
  \mapsto 
  \left[ \begin{smallmatrix} 1 & 0 & 0 \\ 0 & 0 & 1 \\ 0 & 1 & 0 \end{smallmatrix} \right] 
  &\qquad
  &
  \left[ \begin{smallmatrix} 0 & 1 & 0 \\ 0 & 0 & 1 \\ 1 & 0 & 0 \end{smallmatrix} \right] 
  \mapsto
  \left[ \begin{smallmatrix} 0 & 0 & 1 \\ 1 & 0 & 0 \\ 0 & 1 & 0 \end{smallmatrix} \right] 
  &\qquad
  &
  \left[ \begin{smallmatrix} 1 & 0 & 0 \\ 1 & 1 & 0 \\ 0 & 0 & 1 \end{smallmatrix} \right]
  \mapsto
  \left[ \begin{smallmatrix} 0 & 0 & 0 \\ 0 & 1 & 0 \\ 0 & 0 & 0 \end{smallmatrix} \right]
  \\   
  &
  8\colon &\quad
  &
  \left[ \begin{smallmatrix} 0 & 1 & 0 \\ 1 & 0 & 0 \\ 0 & 0 & 1 \end{smallmatrix} \right] 
  \mapsto 
  \left[ \begin{smallmatrix} 0 & 1 & 0 \\ 1 & 0 & 0 \\ 0 & 0 & 1 \end{smallmatrix} \right] 
  &\qquad
  &
  \left[ \begin{smallmatrix} 0 & 1 & 0 \\ 0 & 0 & 1 \\ 1 & 0 & 0 \end{smallmatrix} \right] 
  \mapsto
  \left[ \begin{smallmatrix} 0 & 1 & 0 \\ 0 & 0 & 1 \\ 1 & 0 & 0 \end{smallmatrix} \right] 
  &\qquad
  &
  \left[ \begin{smallmatrix} 1 & 0 & 0 \\ 1 & 1 & 0 \\ 0 & 0 & 1 \end{smallmatrix} \right]
  \mapsto
  \left[ \begin{smallmatrix} 1 & 0 & 0 \\ 0 & 0 & 0 \\ 0 & 0 & 0 \end{smallmatrix} \right]
  \end{alignat*}

\subsection{Representations 9, 10; dimension 6}

In both cases the last two generators are sent to the $6 \times 6$ zero matrix.
The first three generators have the following representations; from now on we use dot for zero
in the representation matrices:
  \begin{alignat*}{4}
  &
  9\colon &\quad
  &
  \left[ \begin{smallmatrix} 0 & 1 & 0 \\ 1 & 0 & 0 \\ 0 & 0 & 1 \end{smallmatrix} \right] 
  \mapsto 
  \left[ \begin{smallmatrix} 
  . & . & 1 & . & . & . \\ 
  . & . & . & . & 1 & . \\ 
  1 & . & . & . & . & . \\ 
  . & . & . & . & . & 1 \\ 
  . & 1 & . & . & . & . \\ 
  . & . & . & 1 & . & . 
  \end{smallmatrix} \right] 
  &\quad
  &
  \left[ \begin{smallmatrix} 0 & 1 & 0 \\ 0 & 0 & 1 \\ 1 & 0 & 0 \end{smallmatrix} \right] 
  \mapsto
  \left[ \begin{smallmatrix} 
  . & . & . & . & 1 & . \\ 
  . & . & 1 & . & . & . \\ 
  . & . & . & . & . & 1 \\ 
  1 & . & . & . & . & . \\ 
  . & . & . & 1 & . & . \\
  . & 1 & . & . & . & .  
  \end{smallmatrix} \right] 
  &\quad
  &
  \left[ \begin{smallmatrix} 1 & 0 & 0 \\ 1 & 1 & 0 \\ 0 & 0 & 1 \end{smallmatrix} \right]
  \mapsto
  \left[ \begin{smallmatrix} 
  1 & . & . & . & . & . \\ 
  . & 1 & . & . & . & . \\ 
  . & . & . & . & . & . \\ 
  . & . & . & 1 & . & . \\ 
  . & . & . & . & . & . \\ 
  . & . & . & . & . & . 
  \end{smallmatrix} \right] 
  \\ 
  &
  10\colon &\quad
  &
  \left[ \begin{smallmatrix} 0 & 1 & 0 \\ 1 & 0 & 0 \\ 0 & 0 & 1 \end{smallmatrix} \right] 
  \mapsto 
  \left[ \begin{smallmatrix} 
  . & . & 1 & . & . & . \\ 
  . & . & . & . & 1 & . \\ 
  1 & . & . & . & . & . \\ 
  . & . & . & . & . & 1 \\ 
  . & 1 & . & . & . & . \\ 
  . & . & . & 1 & . & . 
  \end{smallmatrix} \right] 
  &\quad
  &
  \left[ \begin{smallmatrix} 0 & 1 & 0 \\ 0 & 0 & 1 \\ 1 & 0 & 0 \end{smallmatrix} \right] 
  \mapsto
  \left[ \begin{smallmatrix} 
  . & . & . & . & 1 & . \\ 
  . & . & 1 & . & . & . \\ 
  . & . & . & . & . & 1 \\ 
  1 & . & . & . & . & . \\ 
  . & . & . & 1 & . & . \\
  . & 1 & . & . & . & .  
  \end{smallmatrix} \right] 
  &\quad
  &
  \left[ \begin{smallmatrix} 1 & 0 & 0 \\ 1 & 1 & 0 \\ 0 & 0 & 1 \end{smallmatrix} \right]
  \mapsto
  \left[ \begin{smallmatrix} 
  . & . & . & . & . & . \\ 
  . & 1 & . & . & . & . \\ 
  . & . & . & . & . & . \\ 
  . & . & . & . & . & . \\ 
  . & . & . & . & . & . \\ 
  . & . & . & . & . & . 
  \end{smallmatrix} \right] 
  \end{alignat*} 

\subsection{Representation 11, dimension 7}

This is the smallest representation in which all five generators are sent to nonzero matrices:
  \begin{align*}
  &
  \left[ \begin{smallmatrix} 0 & 1 & 0 \\ 1 & 0 & 0 \\ 0 & 0 & 1 \end{smallmatrix} \right] 
  \mapsto 
  \left[ \begin{smallmatrix} 
  1 & . & . & . & . & . & . \\ 
  . & 1 & . & . & . & . & . \\ 
  . & . & . & . & 1 & . & . \\ 
  . & . & . & . & . & 1 & . \\ 
  . & . & 1 & . & . & . & . \\ 
  . & . & . & 1 & . & . & . \\ 
  . & . & . & . & . & . & 1 
  \end{smallmatrix} \right] 
  \quad
  \left[ \begin{smallmatrix} 0 & 1 & 0 \\ 0 & 0 & 1 \\ 1 & 0 & 0 \end{smallmatrix} \right] 
  \mapsto
  \left[ \begin{smallmatrix} 
  1 & . & . & . & . & . & . \\ 
  . & . & . & . & 1 & . & . \\ 
  . & 1 & . & . & . & . & . \\ 
  . & . & . & . & . & 1 & . \\ 
  . & . & 1 & . & . & . & . \\ 
  . & . & . & . & . & . & 1 \\ 
  . & . & . & 1 & . & . & . 
  \end{smallmatrix} \right] 
  \quad
  \left[ \begin{smallmatrix} 1 & 0 & 0 \\ 1 & 1 & 0 \\ 0 & 0 & 1 \end{smallmatrix} \right]
  \mapsto
  \left[ \begin{smallmatrix} 
  1 & . & . & . & . & . & . \\ 
  . & 1 & . & . & . & . & . \\ 
  . & . & 1 & . & . & . & . \\ 
  . & . & . & 1 & . & . & . \\ 
  . & . & . & . & . & . & . \\ 
  . & . & . & . & . & . & . \\ 
  . & . & . & . & 1 & . & 1 
  \end{smallmatrix} \right] 
  \\ 
  &
  \left[ \begin{smallmatrix} 1 & 0 & 0 \\ 0 & 1 & 0 \\ 0 & 0 & 0 \end{smallmatrix} \right] 
  \mapsto 
  \left[ \begin{smallmatrix} 
  \pp 1 & \pp 1 & \pp . & \pp . & \pp . & \pp . & \pp . \\ 
  \pp . & \pp . & \pp . & \pp . & \pp . & \pp . & \pp . \\ 
  \pp . & \pp . & \pp 1 & \pp 1 & \pp . & \pp . & \pp . \\ 
  \pp . & \pp . & \pp . & \pp . & \pp . & \pp . & \pp . \\ 
  \pp . & \pp . & \pp . & \pp . & \pp 1 & \pp 1 & \pp . \\ 
  \pp . & \pp . & \pp . & \pp . & \pp . & \pp . & \pp . \\ 
  -1 & -1 & -1 & -1 & -1 & -1 & \pp . 
  \end{smallmatrix} \right] 
  \quad
  \left[ \begin{smallmatrix} 0 & 1 & 1 \\ 1 & 0 & 1 \\ 1 & 1 & 0 \end{smallmatrix} \right] 
  \mapsto
  \left[ \begin{smallmatrix} 
  1 & . & . & . & . & . & . \\ 
  . & . & . & . & . & . & . \\ 
  . & . & . & . & . & . & . \\ 
  . & . & . & . & 1 & . & . \\ 
  . & . & . & . & . & . & . \\ 
  . & . & 1 & . & . & . & . \\ 
  . & 1 & . & . & . & . & . 
  \end{smallmatrix} \right] 
  \end{align*} 

\subsection{Representation 12, dimension 9}

In this case and the following two, we omit the generators \ref{fivegenerators} and present only 
the corresponding representation matrices in the order displayed in \eqref{fivegenerators}:
  \begin{align*}
  &
  \left[ \begin{smallmatrix} 
  . & 1 & . &\pp  . &\pp  . &\pp  . &\pp  . & . & . \\ 
  1 & . & . &\pp  . &\pp  . &\pp  . &\pp  . & . & . \\ 
  . & . & 1 &\pp  . &\pp  . &\pp  . &\pp  . & . & . \\ 
  . & . & . &    -1 &\pp  . &\pp  . &\pp  . & . & . \\ 
  . & . & . &\pp  . &\pp  . &    -1 &\pp  . & . & . \\ 
  . & . & . &\pp  . &    -1 &\pp  . &\pp  . & . & . \\ 
  . & . & . &\pp  . &\pp  . &\pp  . &    -1 & . & . \\ 
  . & . & . &\pp  . &\pp  . &\pp  . &\pp  . & . & 1 \\ 
  . & . & . &\pp  . &\pp  . &\pp  . &\pp  . & 1 & . 
  \end{smallmatrix} \right] 
  \qquad
  \left[ \begin{smallmatrix} 
  . &   -1 & . &\pp . &\pp . &\pp . &\pp . &\pp . & . \\ 
  . &\pp . & . &   -1 &\pp . &\pp . &\pp . &\pp . & . \\ 
  . &\pp . & . &\pp . &\pp . &   -1 &\pp . &\pp . & . \\ 
  1 &\pp . & . &\pp . &\pp . &\pp . &\pp . &\pp . & . \\ 
  . &\pp . & 1 &\pp . &\pp . &\pp . &\pp . &\pp . & . \\ 
  . &\pp . & . &\pp . &   -1 &\pp . &\pp . &\pp . & . \\ 
  . &\pp . & . &\pp . &\pp . &\pp . &\pp . &\pp . & 1 \\ 
  . &\pp . & . &\pp . &\pp . &\pp . &   -1 &\pp . & . \\ 
  . &\pp . & . &\pp . &\pp . &\pp . &\pp . &   -1 & . 
  \end{smallmatrix} \right] 
  \\
  &
  \left[ \begin{smallmatrix}  
  1 & . & . & . & . & . & . & . & . \\ 
  . & . & . & . & . & . & . & . & . \\ 
  . & 1 & 1 & . & . & . & . & . & . \\ 
  . & . & . & . & . & . & . & . & . \\ 
  . & . & . & . & . & . & . & . & . \\ 
  . & . & . & . & . & . & . & . & . \\ 
  . & . & . & . & . & . & . & . & . \\ 
  . & . & . & . & . & 1 & . & 1 & . \\ 
  . & . & . & . & . & . & . & . & . 
  \end{smallmatrix} \right] 
  \qquad
  \left[ \begin{smallmatrix} 
  . & . & . & . & . & . & . & . & . \\ 
  . & . & . & . & . & . & . & . & . \\ 
  . & . & . & . & . & . & . & . & . \\ 
  . & . & . & 1 & 1 & 1 & 1 & . & . \\ 
  . & . & . & . & . & . & . & . & . \\ 
  . & . & . & . & . & . & . & . & . \\ 
  . & . & . & . & . & . & . & . & . \\ 
  . & . & . & . & . & . & . & . & . \\ 
  . & . & . & . & . & . & . & . & . 
  \end{smallmatrix} \right] 
  \qquad
  \left[ \begin{smallmatrix} 
  \pp . &\pp . & . &\pp . & . & . & . & . & . \\ 
  \pp . &\pp . & . &\pp . & . & . & . & . & . \\ 
  \pp . &\pp . & . &\pp . & . & . & . & . & . \\ 
  \pp . &\pp . & . &\pp . & . & . & . & . & . \\ 
  \pp . &\pp . & . &\pp . & . & . & . & . & . \\ 
  \pp . &\pp . & . &\pp . & . & . & . & . & . \\ 
  \pp . &\pp . & . &   -1 & . & . & . & . & . \\ 
  \pp . &   -1 & . &\pp . & . & . & . & . & . \\ 
     -1 &\pp . & . &\pp . & . & . & . & . & . 
  \end{smallmatrix} \right] 
  \end{align*} 

\subsection{Representation 13, dimension 9}

  \begin{align*}
  &
  \left[ \begin{smallmatrix} 
  . & 1 & . & . & . & . & . & . & . \\ 
  1 & . & . & . & . & . & . & . & . \\ 
  . & . & 1 & . & . & . & . & . & . \\ 
  . & . & . & 1 & . & . & . & . & . \\ 
  . & . & . & . & . & 1 & . & . & . \\ 
  . & . & . & . & 1 & . & . & . & . \\ 
  . & . & . & . & . & . & 1 & . & . \\ 
  . & . & . & . & . & . & . & . & 1 \\ 
  . & . & . & . & . & . & . & 1 & . 
  \end{smallmatrix} \right] 
  \qquad
  \left[ \begin{smallmatrix} 
  . & 1 & . & . & . & . & . & . & . \\ 
  . & . & . & 1 & . & . & . & . & . \\ 
  . & . & . & . & . & 1 & . & . & . \\ 
  1 & . & . & . & . & . & . & . & . \\ 
  . & . & 1 & . & . & . & . & . & . \\ 
  . & . & . & . & 1 & . & . & . & . \\ 
  . & . & . & . & . & . & . & . & 1 \\ 
  . & . & . & . & . & . & 1 & . & . \\ 
  . & . & . & . & . & . & . & 1 & . 
  \end{smallmatrix} \right] 
  \\
  &
  \left[ \begin{smallmatrix} 
  1 & . & . & . & . & . & . & . & . \\ 
  . & . & . & . & . & . & . & . & . \\ 
  . & 1 & 1 & . & . & . & . & . & . \\ 
  . & . & . & . & . & . & . & . & . \\ 
  . & . & . & . & . & . & . & . & . \\ 
  . & . & . & . & . & . & . & . & . \\ 
  . & . & . & . & . & . & . & . & . \\ 
  . & . & . & . & . & 1 & . & 1 & . \\ 
  . & . & . & . & . & . & . & . & . 
  \end{smallmatrix} \right] 
  \qquad
  \left[ \begin{smallmatrix} 
  . & . & . & . & . & . & . & . & . \\ 
  . & . & . & . & . & . & . & . & . \\ 
  . & . & . & . & . & . & . & . & . \\ 
  . & . & . & 1 & 1 & 1 & 1 & . & . \\ 
  . & . & . & . & . & . & . & . & . \\ 
  . & . & . & . & . & . & . & . & . \\ 
  . & . & . & . & . & . & . & . & . \\ 
  . & . & . & . & . & . & . & . & . \\ 
  . & . & . & . & . & . & . & . & . 
  \end{smallmatrix} \right] 
  \qquad
  \left[ \begin{smallmatrix} 
  . & . & . & . & . & . & . & . & . \\ 
  . & . & . & . & . & . & . & . & . \\ 
  . & . & . & . & . & . & . & . & . \\ 
  . & . & . & . & . & . & . & . & . \\ 
  . & . & . & . & . & . & . & . & . \\ 
  . & . & . & . & . & . & . & . & . \\ 
  . & . & . & 1 & . & . & . & . & . \\ 
  . & 1 & . & . & . & . & . & . & . \\ 
  1 & . & . & . & . & . & . & . & . 
  \end{smallmatrix} \right] 
  \end{align*} 

\subsection{Representation 14, dimension 12}

  \begin{align*}
  &
  \left[ \begin{smallmatrix} 
  1 & . & . & . & . & . & . & . & . & . & . & . \\ 
  . & . & . & . & . & 1 & . & . & . & . & . & . \\ 
  . & . & . & . & . & . & 1 & . & . & . & . & . \\ 
  . & . & . & . & . & . & . & 1 & . & . & . & . \\ 
  . & . & . & . & . & . & . & . & . & . & 1 & . \\ 
  . & 1 & . & . & . & . & . & . & . & . & . & . \\ 
  . & . & 1 & . & . & . & . & . & . & . & . & . \\ 
  . & . & . & 1 & . & . & . & . & . & . & . & . \\ 
  . & . & . & . & . & . & . & . & 1 & . & . & . \\ 
  . & . & . & . & . & . & . & . & . & . & . & 1 \\ 
  . & . & . & . & 1 & . & . & . & . & . & . & . \\ 
  . & . & . & . & . & . & . & . & . & 1 & . & . 
  \end{smallmatrix} \right] 
  \qquad
  \left[ \begin{smallmatrix} 
  . & . & . & . & . & 1 & . & . & . & . & . & . \\ 
  1 & . & . & . & . & . & . & . & . & . & . & . \\ 
  . & . & . & . & . & . & 1 & . & . & . & . & . \\ 
  . & . & . & . & . & . & . & . & . & . & 1 & . \\ 
  . & . & . & . & . & . & . & 1 & . & . & . & . \\ 
  . & 1 & . & . & . & . & . & . & . & . & . & . \\ 
  . & . & . & . & . & . & . & . & 1 & . & . & . \\ 
  . & . & . & . & . & . & . & . & . & . & . & 1 \\ 
  . & . & 1 & . & . & . & . & . & . & . & . & . \\ 
  . & . & . & 1 & . & . & . & . & . & . & . & . \\ 
  . & . & . & . & . & . & . & . & . & 1 & . & . \\ 
  . & . & . & . & 1 & . & . & . & . & . & . & . 
  \end{smallmatrix} \right] 
  \qquad
  \left[ \begin{smallmatrix} 
  1 & . & . & . & . & . & . & . & . & . & . & . \\ 
  . & 1 & . & . & . & . & . & . & . & . & . & . \\ 
  . & . & 1 & . & . & . & . & . & . & . & . & . \\ 
  . & . & . & 1 & . & . & . & . & . & . & . & . \\ 
  . & . & . & . & 1 & . & . & . & . & . & . & . \\ 
  . & . & . & . & . & . & . & . & . & . & . & . \\ 
  . & . & . & . & . & . & . & . & . & . & . & . \\ 
  . & . & . & . & . & . & . & . & . & . & . & . \\ 
  . & . & . & . & . & 1 & . & . & 1 & . & . & . \\ 
  . & . & . & . & . & . & . & . & . & 1 & . & . \\ 
  . & . & . & . & . & . & . & . & . & . & . & . \\ 
  . & . & . & . & . & . & . & . & . & . & . & . 
  \end{smallmatrix} \right]
  \\
  &
  \left[ \begin{smallmatrix} 
  . &\pp . &\pp . &\pp . & . &\pp . &\pp . &\pp . & . & . & . & . \\ 
  . &\pp 1 &\pp 1 &\pp 1 & . &\pp . &\pp . &\pp . & . & . & . & . \\ 
  . &\pp . &\pp . &\pp . & . &\pp . &\pp . &\pp . & . & . & . & . \\ 
  . &\pp . &\pp . &\pp . & . &\pp . &\pp . &\pp . & . & . & . & . \\ 
  . &\pp . &\pp . &\pp . & . &\pp . &\pp . &\pp . & . & . & . & . \\ 
  . &\pp . &\pp . &\pp . & . &\pp 1 &\pp 1 &\pp 1 & . & . & . & . \\ 
  . &\pp . &\pp . &\pp . & . &\pp . &\pp . &\pp . & . & . & . & . \\ 
  . &\pp . &\pp . &\pp . & . &\pp . &\pp . &\pp . & . & . & . & . \\ 
  . &   -1 &   -1 &   -1 & . &   -1 &   -1 &   -1 & . & . & . & . \\ 
  . &\pp 1 &\pp 1 &\pp 1 & . &\pp . &\pp . &\pp . & . & . & . & . \\ 
  . &\pp . &\pp . &\pp . & . &\pp . &\pp . &\pp . & . & . & . & . \\ 
  . &\pp . &\pp . &\pp . & . &\pp 1 &\pp 1 &\pp 1 & . & . & . & . 
  \end{smallmatrix} \right] 
  \qquad
  \left[ \begin{smallmatrix} 
  . & . & . & . & . & . & . & . & . & . & . & . \\ 
  . & . & . & . & . & . & . & . & . & . & . & . \\ 
  . & . & . & . & . & 1 & . & . & . & . & . & . \\ 
  . & . & . & . & . & . & . & . & . & . & . & . \\ 
  . & . & . & . & . & . & . & . & . & . & . & . \\ 
  . & . & . & . & . & . & . & . & . & . & . & . \\ 
  . & 1 & . & . & . & . & . & . & . & . & . & . \\ 
  . & . & . & . & . & . & . & . & . & . & . & . \\ 
  1 & . & . & . & . & . & . & . & . & . & . & . \\ 
  . & . & . & . & . & . & . & . & . & . & . & . \\ 
  . & . & . & . & . & . & . & . & . & . & . & . \\ 
  . & . & . & . & . & . & . & . & . & . & . & . 
  \end{smallmatrix} \right] 
  \end{align*} 

\begin{remark}
Only the representations of dimension $d \ge 7$ are faithful.
It was proved by Kim and Roush \cite{KimRoush} that the smallest faithful representation
of $\mathcal{B}_n$ has dimension $2^n{-}1$; see also \cite{BremnerElBachraoui}.
\end{remark}



\begin{thebibliography}{99}

\bibitem{Bremner}
\textsc{M. R. Bremner}:
How to compute the Wedderburn decomposition of a finite-dimensional associative algebra.
\emph{Groups Complex. Cryptol.} 
3 (2011), no.~1, 47--66. 

\bibitem{BremnerLLL}
\textsc{M. R. Bremner}:
\emph{Lattice Basis Reduction: An Introduction to the LLL Algorithm and Its Applications}. 
Pure and Applied Mathematics, 300. 
CRC Press, Boca Raton, FL, 2012. 
xviii+316 pp. 
ISBN: 978-1-4398-0702-6.

\bibitem{BremnerElBachraoui}
\textsc{M. R. Bremner, M. El Bachraoui}:
On the semigroup algebra of binary relations.
\emph{Comm. Algebra} 
38 (2010), no.~9, 3499--3505.

\bibitem{deCaenGregory}
\textsc{D. de Caen, D. A. Gregory}:
Primes in the semigroup of Boolean matrices.
\emph{Linear Algebra Appl.} 
37 (1981) 119--134. 

\bibitem{Devadze}
\textsc{H. N. Devadze}:
Generating sets of the semigroup of all binary relations in a finite set.
\emph{Dokl. Akad. Nauk BSSR} 
12 (1968) 765--768. 

\bibitem{Dickson}
\textsc{L. E. Dickson}:
\emph{Algebras and Their Arithmetics}. 
Dover Publications, New York, 1960. 
xii+241 pp.
ISBN: 978-0-486-60616-3.
Unaltered republication of first edition (University of Chicago Press, 1923). 

\bibitem{Drazin}
\textsc{M. P. Drazin}:
Maschke's theorem for semigroups.
\emph{J. Algebra} 
72 (1981), no.~1, 269--278. 

\bibitem{IvanyosRonyai}
\textsc{G. Ivanyos, L. R\'onyai}:
Computations in associative and Lie algebras. 
Chapter 5 of \emph{Some Tapas of Computer Algebra}
(A. M. Cohen, H. Cuypers, H. Sterk, editors). 
Springer, 1999.
xiv+352 pp.
ISBN: 978-3-540-63480-5.

\bibitem{IRS}
\textsc{G. Ivanyos, L. R\'onyai, J. Schicho}:
Splitting full matrix algebras over algebraic number fields.
\emph{J. Algebra} 
354 (2012) 211--223. 

\bibitem{KimRoush}
\textsc{K. H. Kim, F. W. Roush}:
Linear representations of semigroups of Boolean matrices.
\emph{Proc. Amer. Math. Soc.} 
63 (1977), no.~2, 203--207.

\bibitem{Konieczny}
\textsc{J. Konieczny}:
A proof of Devadze's theorem on generators of the semigroup of Boolean matrices.
\emph{Semigroup Forum} 
83 (2011), no.~2, 281--288. 

\bibitem{Preston}
\textsc{G. B. Preston}: 
Any group is a maximal subgroup of the semigroup of binary relations on some set.
\emph{Glasg. Math. J.} 
14 (1973) 21--24.

\end{thebibliography}
\end{document}